\documentclass[11pt,leqno,twoside]{amsart}
\usepackage{amsmath}
\usepackage{amsthm}                 
\usepackage{amssymb}

\usepackage{graphicx}

\usepackage{epsfig}
\newtheorem{theorem}{Theorem}
\newtheorem{corollary}{Corollary}
\newtheorem{proposition}{Proposition}

\theoremstyle{definition}

\newtheorem{example}{Example}

\newenvironment{definition}
{\smallskip\noindent{\bf Definition\/}:}{\smallskip\par}

\newenvironment{remark}
{\smallskip\noindent{\bf Remark\/}.}{\smallskip\par}

\usepackage{hhline}

\newcommand{\CC}{{\Bbb C}}

\newcommand{\PP}{{\Bbb P}}

\newcommand{\RR}{{\Bbb R}}
\newcommand{\ZZ}{{\Bbb Z}}

\newcommand{\DD}{{\Bbb D}}

\newcommand{\TT}{{\Bbb T}}
\newcommand{\mathcalO}{{\mathcal O}}

\newcommand{\mathcalJ}{{\mathcal J}}

\newcommand{\eps}{\varepsilon}
\newcommand{\icis}{{ICIS}}
\newcommand{\ind}{{\rm ind}}
\newcommand{\sgn}{{\rm sgn}}
\newcommand{\Sing}{{\rm Sing}}

\newcommand{\cO}{{\mathcal O}}
\newcommand{\bC}{{\mathbb C}}

\begin{document}

\title{Indices of vector fields and 1-forms on singular varieties}
\author{W.~Ebeling}
\address{Universit\"{a}t Hannover, Institut f\"{u}r Algebraische Geometrie, Postfach
6009, D-30060 Hannover, Germany}
 \email{ebeling@math.uni-hannover.de} 
\author{S.~M.~Gusein-Zade}
\address{Moscow State University, Faculty of Mechanics and Mathematics, Moscow,
119992, Russia}
\email{sabir@mccme.ru}
\thanks{This work was partially supported by the DFG-programme ''Global methods in
complex geometry'' (Eb 102/4--3), grants RFBR--04--01--00762,
NWO-RFBR 047.011.2004.026.}

\maketitle

\begin{abstract}
We discuss different generalizations of the classical notion of
the index of a singular point of a vector field to the case of
vector fields or 1-forms on singular varieties, describe relations
between them and formulae for their computation.
\end{abstract}

\section*{Introduction}
An isolated singular point (zero) of a (continuous) vector field
on a smooth manifold has a well known invariant~--- the index.
A neighbourhood of a point $P$ in a manifold $M^n$ of dimension
$n$ can be identified with a neighbourhood of the origin in the
affine (coordinate) space $\RR^n$. A germ $X$ of a vector field
on $(\RR^n, 0)$ can be written as
$X(x) = \sum\limits_{i=1}^n X_i(x)\frac{\partial\ \ }{\partial x_i}$.
Let $B_\eps^{n}(0)$ be the ball of radius $\eps$ centred at
the origin in $\RR^n$ such that the vector field $X$ is defined
in a neighbourhood of this ball and has no zeros in it except at
the origin. Let $S_\eps^{n-1}(0)$ be the $(n-1)$-dimensional sphere
$\partial B_\eps^{n}(0)$. The vector field $X$ defines a map
$$
\frac{X}{\Vert X\Vert}:S_\eps^{n-1}(0)\to S_1^{n-1}\, .
$$
The index ${\ind}_P\,X$ of the vector field $X$ at the point
$P$ is defined as the degree of the map $\frac{X}{\Vert X\Vert}$.
One can see that it is independent of the chosen coordinates.

If the vector field $X$ is smooth and its singular point $P$ is
non-degenerate, i.e. if
$J_{X,P}:= \det\left(\frac{\partial X_i}{\partial x_j}(0)\right) \ne 0$,
then ${\ind}_P\,X= \mbox{sgn\,}J_{X,P}$, i.e. it is equal to $1$
if $J_{X,P}>0$ and is equal to $-1$ if $J_{X,P}<0$. The index of
an arbitrary isolated singular point $P$ of a smooth vector field $X$
is equal to the number of non-degenerate singular points which split
from the point $P$ under a generic deformation $\widetilde X$ of the
vector field $X$ in a neighbourhood of the point $P$ counted with
the corresponding signs ($\mbox{sgn\,}J_{{\widetilde X},{\widetilde P}})$.
This follows, in particular, from the fact that the index of
a vector field satisfies the law of conservation of number (see below).

One of the most important properties of the index of a vector field
is the Poincar\'e--Hopf theorem. Suppose that the manifold $M$ is
closed, i.e. compact without boundary, and that the vector field $X$
has finitely many singular points on it.

\begin{theorem}[Poincar\'e-Hopf]\label{theo1}
The sum
$$
\sum\limits_{P\in {\rm Sing}\, X}{\ind}_P\,X
$$
of indices of singular points of the vector field $X$ is equal to
the Euler characteristic $\chi(M)$ of the manifold $M$.
\end{theorem}

If $M$ is an $n$-dimensional complex analytic (compact) manifold,
then its Euler characteristic $\chi(M)$ is equal to the characteristic
number $\langle c_n(TM), [M]\rangle$, where $c_n(M)$ is the top Chern
class of the manifold $M$. If a vector field $X$ on a complex manifold
$M$ is holomorphic and a singular point $P$ of it is non-degenerate,
then the index ${\ind}_P\,X$ is equal to $+1$. The index of an
isolated singular point $P$ of a holomorphic vector field $X$
is positive. It is equal to the number of non-degenerate singular
points which split from the point $P$ under a generic deformation
of the vector field $X$ in a neighbourhood of the point $P$.

There exists an algebraic formula for the index ${\ind}_P\,X$ of
an isolated singular point of a holomorphic vector field. In local
coordinates centred at the point $P$, let the vector field $X$ be
equal to $\sum\limits_{i=1}^n X_i(x)\frac{\partial\ \ }{\partial x_i}$.
Let ${\mathcal O}_{\CC^n,0}$ be the ring of germs of holomorphic functions
of $n$ variables.
\begin{theorem}\label{theo2}
$$
{\ind}_P X = \dim_\CC {\mathcal O}_{\CC^n,0}/(X_1, \ldots , X_n)\,,
$$
where $(X_1, \ldots , X_n)$ is the ideal generated by the germs
$X_1$, \dots, $X_n$.
\end{theorem}

This statement was known for a long time, however, it seems that
the first complete proof appeared in \cite{Palamodov67}.

An algebraic formula for the index ${\ind}_P\,X$ of an isolated
(in fact for an algebraically isolated: see the definition below)
singular point of a smooth (say, $C^\infty$) vector field was given in
\cite{EL77, Khimshiashvili77}. In local coordinates centred at the point
$P$, let $X=\sum\limits_{i=1}^n X_i(x)\frac{\partial\ \ }{\partial x_i}$.
Let ${\mathcal E}_{\RR^n,0}$ be the ring of germs of smooth ($C^\infty$)
functions of $n$ variables.

\begin{definition}
The singular point $P$ of the vector field $v$ is algebraically
isolated if the factor ring ${\mathcal E}_{\RR^n,0}/(X_1, \ldots , X_n)$
has a finite dimension as a vector space.
\end{definition}

If a singular point of a vector field is algebraically isolated, then
it is isolated. If, in the local coordinates $x_1$, \dots, $x_n$, the
vector field $X$ is (real) analytic (i.e., if the components $X_i$ are
analytic functions of $x_1$, \dots, $x_n$), one can consider its
complexification $X_\CC$ which is a holomorphic vector field in a
neighbourhood of the origin in $\CC^n$. In this case the singular point
of the vector field $X$ is algebraically isolated if and only if the
singular point of the vector field $X_\CC$ (the origin) is isolated. Let
$$
J_X=J_X(x):=\det\left(\frac{\partial X_i}{\partial x_j}\right)
$$
be the {\it Jacobian} of the vector field $X$. One can prove that
$J_X\ne 0$ in the local ring $R_X={\mathcal E}_{\RR^n,0}/(X_1,\ldots,X_n)$.
Moreover, the ring $R_X$
has a one dimensional minimal ideal and this ideal is generated by the
Jacobian $J_X$. Let $\ell:R_X\to \RR$ be a linear function on the
ring $R_X$ (considered as a vector space) such that $\ell(J_X)>0$.
Consider the quadratic form $Q$ on $R_X$ defined by
$$
Q(\varphi, \psi)=\ell(\varphi\cdot\psi)\,.
$$

\begin{theorem}[Eisenbud-Levine-Khimshiashvili]\label{theo3}
The index ${\ind}_P X$ of the singular point $P$ of the vector field $X$
is equal to the signature $\sgn \, Q$ of the quadratic form $Q$.
\end{theorem}

For a proof of this theorem see also \cite{AGV85}.

In what follows we shall, in particular, discuss generalizations of
the notion of the index of a vector field to real or complex analytic
variety with singularities and problems of their computation.


\section{1-forms versus vector fields}\label{sec1}
Traditionally the definition of the index and the corresponding results
(say, the Poincar\'e--Hopf and the Eisenbud--Levine--Khimshiashvili
theorems) are formulated for vector fields. However, instead of vector
fields one can consider 1-forms. Using a Riemannian metric one can identify
vector fields and 1-forms on a smooth ($C^\infty$) manifold. Therefore on
a smooth manifold all notions and statements concerning vector fields
can be formulated for 1-forms as well. In particular, the index
${\ind}_P\,\omega$ of a 1-form $\omega$ on a smooth $n$-dimensional
manifold $M$ at an isolated singular point (zero) $P$ is defined (in
local coordinates) as the degree of the map
$\frac{\omega}{\Vert \omega\Vert}:S_\eps^{n-1}\to S_1^{n-1}$,
where $S_1^{n-1}$ is the unit sphere in the dual space. In the local
complex analytic
situation one can identify germs of complex analytic vector fields and
germs of complex analytic 1-forms using local coordinates: the vector
field $\sum\limits X_i\frac{\partial\ }{\partial z_i}$ corresponds to
the 1-form $\sum\limits X_i dz_i$. In the Poincar\'e--Hopf theorem
for 1-forms on an $n$-dimensional complex manifold $M$ the sum of
indices of singular points of a 1-form is equal to
$\langle c_n(T^*M), [M]\rangle= (-1)^n\chi(M)$. The only essential
difference between vector fields and 1-forms in the
smooth complex analytic case is that
non-trivial complex analytic (global) vector fields and 1-forms exist,
generally speaking, on different complex analytic manifolds: in the
odd-dimensional case vector fields (respectively 1-forms) can exist
only on manifolds with non-negative (respectively non-positive) Euler
characteristic.

However, for the case of singular varieties the situation becomes quite
different. In a series of papers instead of considering vector fields
on a variety we started to consider 1-forms. For 1-forms on germs of
singular varieties (real or complex analytic) one can define several
notions of indices usually corresponding to appropriate analogues for
vector fields. However, the properties of the indices for vector fields
and for 1-forms are different. For example, just as above, to an analytic
vector field $X=\sum\limits_{i=1}^N X_i\frac{\partial\ }{\partial x_i}$
on a real variety $(V, 0)\subset (\RR^N, 0)$ (or on a complex analytic
variety $(V, 0)\subset (\CC^N, 0)$) one can associate the 1-form
$\omega=\sum\limits_{i=1}^N X_i dx_i$ (dependent on the choice of
coordinates $x_1$, \dots, $x_N$ in $(\RR^N, 0)$ or in $(\CC^N, 0)$).
If the vector field $X$ has an isolated singular point on $(V,0)$, then,
for a generic (!) choice of the coordinates $x_1$, \dots, $x_N$, the
corresponding 1-form $\omega$ has an isolated singular point as well.
This correspondence does not work in the other direction. Moreover,
whereas for the radial index on a real analytic variety (see below) the
index of a vector field coincides with that of the corresponding 1-form,
this does not hold for other indices, say, for the so called GSV index
(on an isolated complete intersection singularity: {\icis}). In some
cases one can say that the notion of an index of a 1-form (say, of
a holomorphic one) is somewhat more natural than that of a vector field.
For example, the notion of the GSV index of a 1-form (see below) is
``more complex analytic'' (does not use the complex conjugation for the
definition) and ``more geometric'' (uses only objects of the same tensor
type). Moreover the index of an isolated singular point of a holomorphic
1-form on a complex {\icis} can be described as the dimension of an
appropriate algebra. Finally, on a real {\icis} the real index of a 1-form
which is the differential of a germ of a function with an algebraically
isolated singular point plus-minus the Euler characteristic of a (real)
smoothing of the {\icis} can be expressed in terms of the signature of a
quadratic form on a space (certain space of thimbles), the dimension of
which is equal to the (complex) index of the corresponding complexification.

The condition for a vector field to be tangent to a germ of a singular
variety is a very restrictive one. For example, holomorphic vector
fields with isolated zeros exist on a complex analytic variety with
an isolated singularity \cite{BG94}, but not in general.

\begin{example}
Consider the surface $X$ in $\CC^3$ given by the equation
$$
xy(x-y)(x+zy)=0.
$$
It has singularities on the line $\{x=y=0\}$. It
can be considered as a family of four lines in the $(x,y)$-plane with
different cross ratios. Then any holomorphic vector field tangent to
$X$ vanishes on the line $\{x=y=0\}$, because the translation along
such a vector field has to preserve the cross ratio of the lines.
\end{example}

However, 1-forms with isolated singular points always exist.

The idea to consider indices of 1-forms instead of indices of
vector fields (in some situations) was first formulated by
V.I.Arnold (\cite{Arnold79}, see also \cite{Arnold04}).


\section{Radial index} \label{sec2} 
One can say that the notion of the radial index of a vector field or
of a 1-form on a singular variety is a straightforward generalization
of the usual index inspired by the Poincar\'e--Hopf theorem.

Let us start with the setting of (real: $C^\infty$) manifolds with
isolated singularities.

A {\em manifold with isolated singularities} is a topological
space $M$ which has the structure of a smooth (say, $C^\infty$--) manifold
outside of a discrete set $S$ (the {\em set of singular points} of $M$).
A {\em diffeomorphism} between two such manifolds is a homeomorphism
which sends the set of singular points onto the set of singular points and
is a diffeomorphism outside of them. We say that $M$ has a {\em
cone-like singularity} at a (singular) point $P\in S$ if there exists a
neighbourhood of the point $P$ diffeomorphic to the cone over a smooth
manifold $W_P$ ($W_P$ is called the {\em link} of the point $P$). In
what follows we assume all manifolds to have only cone-like
singularities. A (smooth or continuous) {\em vector field} or {\em 1-form}
on a manifold $M$ with isolated singularities is a (smooth or continuous)
vector field or 1-form respectively on the set $M\setminus S$ of regular
points of $M$. The {\em set of singular points} $S_X$ of a vector field $X$
or the {\em set of singular points} $S_\omega$ of a 1-form $\omega$
on a (singular) manifold $M$ is the union of the set of usual singular
points of $X$ or of $\omega$ respectively on $M\setminus S$ (i.e., points at
which $X$ or $\omega$ respectively tends to zero) and of the set $S$ of
singular points of $M$ itself.

For an isolated {\em usual} singular point $P$ of a vector field $X$
or of a 1-form $\omega$ there is defined its index ${\ind}_P\,X$
or ${\ind}_P\, \omega$ respectively.
If the manifold $M$ is closed and has no singularities ($S=\emptyset$)
and the vector field $X$ or the 1-form $\omega$ on $M$ has only isolated
singularities, then
the sum of these indices over all singular points is equal to
the Euler characteristic $\chi(M)$ of the manifold $M$.

Let $(M, P)$ be a cone-like singularity (i.e., a germ of a manifold with
such a singular point). Let $X$ be a vector field defined on an open
neighbourhood $U$ of the point $P$. Suppose that $X$ has no singular
points on $U\setminus\{P\}$. Let $V$ be a closed cone--like neighbourhood
of the point $P$ in $U$ ($V\cong CW_P$, $V\subset U$). On the cone
$CW_P=(W_P \times I)/(W_P \times \{0\})$ ($I=[0, 1]$), there is defined
a natural vector field $\partial/\partial t$ ($t$ is the coordinate on
$I$). Let $X_{\rm rad}$ be the corresponding vector field on $V$. Let
$\widetilde X$ be a
(continuous) vector field on $U$ which coincides with $X$ near
the boundary $\partial U$ of the neighbourhood $U$ and with $X_{\rm rad}$
on $V$ and has only isolated singular points.

\begin{definition}
The {\em radial index} ${\ind}_{\rm rad}\, (X;M,P)$ (or simply
${\ind}_{\rm rad}\,X$) of the vector field $X$ at the point $P$ is equal to
$$
1+\sum_{{\widetilde P}\in S_{\widetilde X}\setminus\{P\}}
{\ind}_{\widetilde P}\widetilde X
$$
(the sum is over all singular points ${\widetilde P}$ of the vector
field $\widetilde X$ except $P$ itself).
\end{definition}

Analogously, we define the radial index of a 1-form $\omega$
defined on an open neighbourhood $U$ of the point $P$ which has
no singular points on $U\setminus\{P\}$. Let $V$ again be a closed
cone--like neighbourhood of $P$ in $U$ ($V\cong CW_P$, $V\subset U$). On
the cone $CW_P=(W_P\times I) /(W_P\times\{0\})$ ($I=[0, 1]$) there is
defined a natural 1-form $dt$. Let $\omega_{\rm rad}$ be the corresponding
1-form on $V$. Let $\widetilde\omega$ be a
(continuous) 1-form on $U$ which
coincides with $\omega$ near the boundary $\partial U$ of the neighbourhood
$U$ and with $\omega_{\rm rad}$ on $V$ and has only isolated singular points.

\begin{definition}
The {\em radial index} ${\ind}_{\rm rad}\, (\omega;M,P)$ (or simply
${\ind}_{\rm rad}\, \omega$) of the 1-form $\omega$ at the point $P$
is equal to
$$
1+\sum_{{\widetilde P}\in S_{\widetilde\omega}\setminus\{P\}}
{\ind}_{\widetilde P}\,\widetilde\omega
$$
(the sum is over all singular points ${\widetilde P}$ of the 1-form
$\widetilde\omega$ except $P$ itself).
\end{definition}

For a cone-like singularity at a point $P\in S$, the link $W_P$ and thus
the cone structure of a neighbourhood are, generally speaking, not
well-defined (cones over different manifolds may be {\em locally}
diffeomorphic). However it is not difficult to show that the indices
${\ind}_{\rm rad}\, X$ and ${\ind}_{\rm rad}\, \omega$ do not depend on
the choice of a cone structure on a neighbourhood and on the choice of
the vector field $\widetilde X$ and the 1-form $\widetilde\omega$
respectively.

\begin{example}
The index of the radial vector field $X_{\rm rad}$ and of the radial
1-form $\omega_{\rm rad}$ is equal to $1$. The index of the vector field
$(-X_{\rm rad})$ and of the 1-form $(-\omega_{\rm rad})$ is equal to
$1-\chi(W_P)$ where $W_P$ is the link of the singular point $P$.
\end{example}

\begin{proposition}\label{prop1}
For a vector field $X$ or a 1-form $\omega$ with isolated singular points
on a closed manifold $M$ with isolated singularities,
the statement of Theorem~\ref{theo1} holds.
\end{proposition}

Now let $(V, 0)\subset(\RR^N, 0)$ be the germ of a real analytic variety,
generally speaking, with a non-isolated singular point at the origin.
Let $V = \bigcup_{i=1}^q V_i$ be a Whitney stratification of the
germ $(V, 0)$.
Let $X$ be a continuous vector field on $(V,0)$ (i.e., the
restriction of a continuous vector field on $(\RR^N, 0)$ tangent to
$V$ at each point) which has an isolated zero at the origin (on $V$).
(Tangency to the stratified set $V$ means that for each point $Q\in V$
the vector $X(Q)$ is tangent to the stratum $V_i$ which contains the point
$Q$.) Let $V_i$ be a stratum of dimension $k$ and let $Q$ be a point of
$V_i$. A neighbourhood of the point $Q$ in $V$ is diffeomorphic to the
direct product of a linear space $\RR^k$ and the cone $CW_Q$ over a
compact singular analytic variety $W_Q$. (Here a diffeomorphism between two
stratified spaces is a homeomorphism which is a diffeomorphism on each
stratum.)
Let $\eps>0$ be small enough so that
in the closed ball $B_\eps$ of radius $\eps$ centred at the origin in
$\RR^N$ the vector field $X$ has no zeros on $V\setminus\{0\}$.
It is
not difficult to show that there exists a (continuous) vector field
$\widetilde X$ on $(V, 0)$ such that:
\begin{enumerate}
\item The vector field $\widetilde X$ is defined on $V\cap B_\eps$.
\item $\widetilde X$ coincides with the vector field $X$ in a
neighbourhood of $V\cap S_\eps$ in $V\cap B_\eps$ ($S_\eps$ is
the sphere $\partial B_\eps$).
\item The vector field $\widetilde X$ has only a finite number of zeros.
\item Each point $Q\in U(0)$ with $\widetilde X(Q)=0$ has a neighbourhood
diffeomorphic to $(\RR^k,0)\times CW_Q$ in which $\widetilde X(y,z)$
($y\in \RR^k$, $z\in CW_Q$) is of the form $Y(y)+X_{\rm rad}(z)$, where
$Y$ is a germ of a vector field on $(\RR^k,0)$ with an isolated singular
point at the origin, $X_{\rm rad}$ is the radial vector field on the
cone $CW_Q$.
\end{enumerate}

Let $S_{\widetilde X}$ be the set of zeros of the vector field
$\widetilde X$ on $V\cap B_\eps$.
For a point $Q\in S_{\widetilde X}$,
let $\widetilde{\ind}(Q):= {\ind}\, (Y;\RR^k,0)$, where $Y$ is
the vector field on $(\RR^k, 0)$ described above.
If $k=0$ (this happens at the origin if it is a stratum of the
stratification $\{V_i\}$), we set $\widetilde{\ind}(Q):= 1$.

\begin{definition} The {\em radial index} ${\ind}_{\rm rad}\, (X;V,0)$
of the vector field $X$ on the variety $V$ at the origin is the number
$$
\sum\limits_{Q\in S_{\widetilde X}} \widetilde{\ind}(Q).
$$
\end{definition}

One can say that the idea of this definition goes back to M.-H.~Schwartz
who defined this index for so called radial vector fields
\cite{Schwartz65, Schwartz86a, Schwartz86b}; cf. also \cite{KT99}.
It was used to define characteristic classes for singular varieties, see e.g. the surveys
\cite{Brasselet00, Seade02}.

Let $\omega$ be an (arbitrary continuous) 1-form
on a neighbourhood of the origin in $\RR^N$ with an isolated singular
point on $(V, 0)$ at the origin.
(A point $Q\in V$ is a singular point of the 1-form $\omega$ on $V$
if it is singular for the restriction of the 1-form $\omega$ to the
stratum $V_i$ which contains the point $Q$; if the stratum $V_i$
is zero-dimensional, the point $Q$ is always singular.)
Let $\eps>0$ be small enough so that
in the closed ball $B_\eps$ of radius $\eps$ centred at the origin in
$\RR^N$ the 1-form $\omega$ has no singular points on $V\setminus\{0\}$.
We say that a 1-form $\omega$ on a germ $(W,0)$ is radial if for any
analytic curve $\gamma:(\RR, 0)\to (W,0)$ different from the trivial one
($\gamma(t)\not\equiv 0$) the value $\omega(\dot\gamma(t))$ of the
1-form $\omega$ on the tangent vector to the curve $\gamma$ is positive
for positive $t$ small enough.
It is easy to see that there exists a 1-form $\widetilde\omega$ on $\RR^N$
such that:
\begin{enumerate}
\item The 1-form $\widetilde\omega$ coincides with the 1-form $\omega$
on a neighbourhood of the sphere $S_\eps = \partial B_\eps$.
\item the vector field $\widetilde X$ has only a finite number of zeros;
\item In a neighbourhood of each singular point
$Q\in (V\cap B_\eps)\setminus\{0\}$, $Q \in V_i$, $\dim V_i=k$,
the 1-form $\widetilde{\omega}$ looks as follows. There exists a (local)
analytic diffeomorphism $h: (\RR^N, \RR^k,0) \to (\RR^N,V_i, Q)$ such that
$h^\ast \widetilde{\omega} =
\pi_1^\ast \widetilde{\omega}_1 + \pi_2^\ast \widetilde{\omega}_2$,
where $\pi_1$ and $\pi_2$ are the natural projections
$\pi_1: \RR^N \to \RR^k$ and $\pi_2: \RR^N \to \RR^{N-k}$ respectively,
$\widetilde{\omega}_1$ is the germ of a 1-form on $(\RR^k,0)$ with an
isolated singular point at the origin, and $\widetilde{\omega}_2$ is
a radial 1-form on $(\RR^{N-k},0)$.
\end{enumerate}

\begin{remark} One can demand that the 1-form $\widetilde\omega_1$ has
a non-degenerate singular point (and therefore
${\ind}\, (\widetilde\omega_1,\RR^k,0) = \pm 1$),
however, this is not necessary for the definition. \end{remark}

Let $S_{\widetilde\omega}$ be the set of singular points of the
1-form $\widetilde\omega$ on $V\cap B_\eps$. For a point
$Q\in S_{\widetilde\omega}$, let
$\widetilde{\ind}(Q):= {\ind}\, (\widetilde\omega_1;\RR^k,0)$.
If $k=0$, we set $\widetilde{\ind}(Q):= 1$.


\begin{definition} The {\em radial index} ${\ind}_{\rm rad}\,
(\omega;V,0)$ of the 1-form $\omega$ on the variety $V$ at the origin is
the sum
$$
\sum\limits_{Q\in S_{\widetilde\omega}} \widetilde{\ind}(Q).
$$
\end{definition}

One can show that these notions are well defined (for the radial index of a
1-form see \cite{EGGD}).

Just because of the definition, the (radial) index satisfies the law
of conservation of number. For a 1-form this means the following: if a
1-form $\omega'$ with isolated singular points on $V$ is close to the
1-form $\omega$, then
$${\ind}_{\rm rad}\,(\omega;V,0) = \sum_{Q \in {\rm
Sing}\, \omega'} {\ind}_{\rm rad}\, (\omega';V,0)$$
where the sum on the
right hand side is over all singular points $Q$ of the 1-form $\omega'$ on
$V$ in a neighbourhood of the origin.

The radial index generalizes the usual index for vector fields or 1-forms
with isolated singularities on a smooth manifold. In particular one has a
generalization of the Poincar\'e-Hopf theorem:

\begin{theorem}[Poincar\'e-Hopf]
For a compact real analytic variety $V$ and a vector field $X$ or
a 1-form $\omega$ with isolated singular points on $V$, one has
$$
\sum_Q {\ind}_{\rm rad}\,(X;V,Q) =
\sum_Q {\ind}_{\rm rad}\,(\omega;V,Q) = \chi(V)
$$
where $\chi(V)$ denotes the Euler characteristic of the set (variety) $V$.
\end{theorem}

In the case of radial vector fields, this theorem is due to M.-H.~Schwartz
\cite{Schwartz65, Schwartz86a, Schwartz86b, Schwartz91} (see also
\cite{BS81, ASV98, KT99}). For a proof of this theorem for the case of
1-forms see \cite{EGGD}.

Now let $(V,0) \subset (\CC^N,0)$ be the germ of a complex analytic
variety of pure dimension $n$.
Let $\omega$ be a (complex and,
generally speaking, continuous) 1-form on a neighbourhood of the
origin in $\CC^N$. In fact there is a one-to-one correspondence between
complex 1-forms on a complex manifold $M^n$ (say, on $\CC^N$) and real
1-forms on it (considered as a real $2n$-dimensional manifold). Namely,
to a complex 1-form $\omega$ one associates the real 1-form
$\eta= {\rm Re}\, \omega$; the 1-form $\omega$ can be restored from
$\eta$ by the formula $\omega(v)=\eta(v)-i\eta(iv)$ for $v \in T_x M^n$.
This means that the index of the real 1-form ${\rm Re}\, \omega$ is an
invariant of the complex 1-form $\omega$ itself. However, on a smooth
manifold ${\ind}_{M^n,x}\, {\rm Re}\, \omega$ does not coincide
with the usual index of the singular point $x$ of the 1-form $\omega$,
but differs from it by the coefficient $(-1)^n$. (E.g., the index of the
(complex analytic) 1-form
$\omega = \sum_{j=1}^n x_j dx_j$ ($(x_1,\ldots, x_n)$ being the
coordinates of $\CC^n$) is equal to 1, whence the index of the real 1-form
${\rm Re}\, \omega = \sum_{j=1}^n u_j du_j - \sum_{j=1}^n v_j dv_j$
($x_j=u_j+iv_j$) is equal to $(-1)^n$.) This explains the following
definition.

\begin{definition}
The {\em (complex radial) index} ${\ind}_{\rm rad}^{\CC}\,(\omega;V,0)$
of the complex 1-form $\omega$ on an $n$-dimensional variety $V$ at the
origin is $(-1)^n$ times the index of the real 1-form ${\rm Re}\,\omega$
on $V$:
$$
{\ind}_{\rm rad}^{\CC} \, (\omega;V,0) =
(-1)^n \,{\ind}_{\rm rad} \, ({\rm Re}\, \omega;V,0).
$$
\end{definition}


\section{GSV index} \label{sec3} 
This generalization of the index makes sense for varieties which have
isolated complete intersection singularities (in particular, for
hypersurfaces with isolated singularities).

Let $V\subset(\CC^{n+k}, 0)$ be an $n$-dimensional isolated complete
intersection singularity ({\icis}) defined by equations $f_1=\ldots=f_k=0$
(i.e.\ $V=f^{-1}(0)$, where $f=(f_1, \ldots, f_k): (\CC^{n+k}, 0)\to
(\CC^k, 0)$, $\dim V=n$) and let $X=\sum\limits_{i=1}^{n+k} X_i
\frac{\partial}{\partial z_i}$ be a germ of
a (continuous) vector field on $(\CC^{n+k}, 0)$ tangent to the {\icis}
$V$. (The latter means that $X(z)\in T_z V$ for all $z\in V\setminus\{0\}$.)
Suppose that the vector field $X$ does not vanish on $V$ in a punctured
neighbourhood of the origin. In this situation the following index
(called the GSV index after X.~G\'omez-Mont, J.~Seade, and A.~Verjovsky)
is defined.

Let $B_\eps\subset \CC^{n+k}$ be the ball of radius $\eps$ centred at
the origin with (positive) $\eps$ small enough so that all the functions
$f_i$ ($i=1,\ldots, k$) and the vector field $X$ are defined in a
neighbourhood of $B_\eps$, $V$ is transversal to the sphere
$S_\eps=\partial B_\eps$, and the vector field $X$ has no zeros on $V$
inside the ball $B_\eps$ except (possibly) at the origin. Let
$K = V \cap S_\eps$ be the link of the {\icis} $(V,0)$. The
manifold $K$ is $(2n-1)$-dimensional and has a natural orientation as
the boundary of the complex manifold $V \cap B_\eps\setminus\{0\}$.

Let ${\mathcal M}(p,q)$, $p\ge q$, be the space of $p \times q$ matrices
with complex
entries and let $D_{p,q}$ be the subspace of ${\mathcal M}(p,q)$ consisting
of matrices of rank less than $q$. The subset $D_{p,q}$ is an
irreducible subvariety of ${\mathcal M}(p,q)$ of codimension $p-q+1$.
The complement $W_{p,q} = {\mathcal M}(p,q) \setminus D_{p,q}$ is the
Stiefel manifold of $q$-frames (collections of $q$ linearly independent
vectors) in $\CC^p$. It is known that the Stiefel manifold $W_{p,q}$
is $2(p-q)$-connected and $H_{2(p-q)+1}(W_{p,q})\cong \ZZ$ (see, e.g.,
\cite{Husemoller75}). The latter fact also proves that $D_{p,q}$ is
irreducible. Since $W_{p,q}$ is the complement of an irreducible complex
analytic subvariety of codimension $p-q+1$ in ${\mathcal M}(p,q)\cong\CC^{pq}$,
there is a natural choice of a generator of the homology group
$H_{2(p-q)+1}(W_{p,q}) \cong \ZZ$. Namely, the (''positive'')
generator is the boundary of a small ball in a smooth complex
analytic slice transversal to $D_{p,q}$ at a non-singular point.
Therefore a map from an oriented smooth ($C^\infty$) closed
$(2(p-q)+1)$-dimensional manifold to $W_{p, q}$ has a degree.

Define the {\em gradient vector field} ${\mbox{grad}}\, f_i$ of a
function germ $f_i$ by
$$
{\mbox{grad}}\, f_i =
\left(\,\overline{\frac{\partial f_i}{\partial z_1}}, \ldots ,
\overline{\frac{\partial f_i}{\partial z_{n+k}}}\,\right)
$$
(${\mbox{grad}}\, f_i$ depends on the choice of the coordinates
$z_1$, \dots, $z_{n+k}$). One has a map
\begin{equation}\label{eq_psi}
\Psi=(X,{\mbox{grad\,}} f_1, \ldots , {\mbox{grad\,}} f_k):
K\to W_{n+k, k+1}
\end{equation}
from the link $K$ to the Stiefel manifold $W_{n+k, k+1}$.

\begin{definition}
The {\em GSV index } $\ind_{\rm GSV} \, (X;V,0)$ of the vector field $X$
on the {\icis} $V$ at the origin is the degree of the map
$$
\Psi: K\to W_{n+k,k+1}\,.
$$
\end{definition}

\begin{remark}
Note that one uses the complex conjugation for this definition
and the components of the discussed map are of different tensor nature.
Whereas $X$ is a vector field, ${\mbox{grad\,}}f_i$ is more similar
to a covector.
\end{remark}

This index was first defined in \cite{GSV91} for vector fields on
isolated hypersurface singularities. In \cite{SS96} it was generalized
to vector fields on {\icis}.

It is convenient to consider the map $\Psi$ as a map from $V$ to
${\mathcal M}(n+k, k+1)$ defined by the formula (\ref{eq_psi}) (in a
neighbourhood of the ball $B_\eps$). It maps the complement of the
origin in $V$ to the Stiefel manifold $W_{n+k, k+1}$. This description
leads to the following definition of the GSV index.

\begin{proposition}\label{prop_gsv1}
The GSV index $\ind_{\rm GSV} \, (X;V,0)$ of the vector field $X$
on the {\icis} $V$ at the origin is equal to the intersection number
$(\Psi(V)\circ D_{n+k, k+1})$ of the image $\Psi(V)$ of the {\icis}
$V$ under the map $\Psi$ and the variety $D_{n+k, k+1}$ at the origin.
\end{proposition}

Note that, even if the vector field $X$ is holomorphic,
$\Psi(V)$ is not, generally speaking, a complex analytic variety.

Another definition/description of the GSV index $\ind_{\rm GSV} \, (X;V,0)$ can
be given in the following way. Let $\widetilde V=V_t=f^{-1}(t)\cap B_\eps$,
where $t\in\CC^k$, $0<\Vert t\Vert\ll\eps$, be the Milnor fibre of the
{\icis} $V$, i.e. a smoothing of it. For $t$ small enough in a neighbourhood
of the sphere $S_\eps=\partial B_\eps$ the manifolds $V=V_0$ and $V_t$
"almost coincide". This gives an up to isotopy well defined  vector field
$\widetilde X$ on the manifold $\widetilde V$ in a neighbourhood of its
boundary $f^{-1}(t)\cap S_\eps$ (generally speaking, not a complex
analytic one). Let us extend this vector field to a vector field (also
denoted by $\widetilde X$) on the entire manifold $\widetilde V$ with only
isolated zeros.

\begin{proposition}\label{prop_gsv2}
One has
$$\ind_{\rm GSV} \, (X;V,0)=\sum\limits_{Q\in \Sing\,\widetilde X}
{\ind}_{Q}\widetilde X.$$
\end{proposition}

This definition can be easily generalized to a germ of a complex analytic
variety with an isolated singularity and with a fixed smoothing.
For example, J.~Seade defined in this way an index for a singular
point of a vector field on a complex analytic surface with a normal
smoothable Gorenstein singularity \cite{Seade87}.
For a singularity which is not an {\icis}, however, it is
possible that such a
smoothing does not exist or there may be a number of different
smoothings what leads to the situation that the index is not well
defined. These difficulties cannot be met for curve singularities.
Thus in this situation the corresponding index is well defined
(cf. \cite{Goryunov00}).
For a general variety with an isolated
singularity, this ambiguity can be avoided for germs of functions
(i.e. for their differentials in our terms) by defining the corresponding
index as a certain residue, what is done in \cite{IS03}.

Now let $V$ be a compact (say, projective) variety all singular points
of which are local {\icis} and let $X$ be a vector field on $V$ with
isolated singular points. One has the following statement.

\begin{proposition}\label{prop_gsv3}
One has
$$\sum\limits_{Q\in \Sing\,\omega}{\ind}_{\rm GSV}(\omega; V, Q)=
\chi(\widetilde V),$$
where $\widetilde V$ is a smoothing of the
variety (local complete intersection) $V$.
\end{proposition}

A similar construction can be considered in the real setting. Namely,
let $V=f^{-1}(0)\subset(\RR^{n+k}, 0)$ be a germ of a real $n$-dimensional
{\icis}, $f=(f_1, \ldots, f_k): (\RR^{n+k}, 0) \to (\RR^k, 0)$ is a real
analytic map, and let $X=\sum\limits_{i=1}^{n+k} X_i dz_i$ be a germ of
a (continuous) vector field on $(\RR^{n+k}, 0)$ tangent to the variety $V$
outside of the origin. Just as above one defines a map $\Psi$ from
the link $K = V \cap S^{n+k-1}_\eps$ of the {\icis} $(V,0)$ to the
Stiefel manifold $W_{n+k,k+1}^\RR$ of $(k+1)$-frames in $\RR^{n+k}$.
The Stiefel manifold $W_{n+k,k+1}^\RR$ is $(n-2)$-connected and its
first non-trivial homology group $H_{n-1}(W_{n+k,k+1}^\RR;\ZZ)$ is isomorphic
to the group $\ZZ$ of integers for $n$ odd and to the group $\ZZ_2$
of order $2$ for $n$ even. Therefore the map $\Psi$ has a degree
defined as an integer or as an integer modulo $2$ depending on the
parity of the dimension $n$. If the manifold $V\setminus\{0\}$ is
not connected, the construction can be applied to each connected
component of it giving a set of degrees (a multi-degree). This
invariant (called {\em real GSV index} of the vector field $X$ on the real
{\icis} $V$) was introduced and studied in \cite{ASV98}.

In \cite{EG01, EG03a} the notion of the GSV index was adapted to
the case of a 1-form. Let $\omega=\sum A_idx_i$ ($A_i=A_i(x)$) be
a germ of a continuous 1-form on $(\CC^{n+k}, 0)$ which as a 1-form
on the {\icis} $V$ has (at most) an isolated singular point at the origin
(thus it does not vanish on the tangent space $T_PV$ to the variety
$V$ at all points $P$ from a punctured neighbourhood of the origin
in $V$). The 1-forms $\omega$, $df_1$, \dots , $df_k$ are linearly
independent for all $P \in K$. Thus one has a map
$$
\Psi=(\omega, df_1, \ldots , df_k) : K \to W_{n+k, k+1}.
$$
Here $W_{n+k, k+1}$ is the Stiefel manifold of $(k+1)$-frames in
the space dual to $\CC^{n+k}$.

\begin{definition}
The {\em GSV index} $\ind_{\rm GSV} \, (\omega;V,0)$ of the 1-form $\omega$
on the {\icis} $V$ at the origin is the degree of the map
$$
\Psi: K\to W_{n+k,k+1}\,.
$$
\end{definition}

Just as above $\Psi$ can be considered as a map from the {\icis} $V$
to the space ${\mathcal M}(n+k, k+1)$ of $(n+k)\times (k+1)$ matrices.
If the 1-form $\omega$ is holomorphic, the map $\Psi$ and the set
$\Psi(V)$ are complex analytic.

The obvious analogues of Propositions~\ref{prop_gsv1}, \ref{prop_gsv2},
and \ref{prop_gsv3} hold.

There exists an algebraic formula for the index
$\ind_{\rm GSV} \, (\omega;V,0)$ of a holomorphic 1-form $\omega$ on
an {\icis} $V$ which gives it as the dimension of a certain algebra
(see Section~\ref{sec6}).

In \cite{BSS05a}, there was defined a generalization of the notion
of the GSV index for a 1-form on a complete intersection singularity
$V=f^{-1}(0)\subset(\CC^{n+k}, 0)$ if the map
$f=(f_1, \ldots, f_k): (\CC^{n+k}, 0) \to (\CC^k, 0)$, defining the
singularity satisfies Thom's $a_f$ regularity condition. (This
condition implies that the Milnor fibre of the complete intersection
singularity $V$ is well-defined.)


\section{Homological index}\label{sec4}
Let $(V,0)\subset (\CC^N, 0)$ be any germ of an analytic variety of pure
dimension $n$ with an isolated singular point at the origin. Suppose $X$
is a complex analytic
vector field tangent to $(V,0)$ with an isolated singular point at
the origin.

Let $\Omega^k_{V,0}$ be be the module of germs of differentiable
$k$-forms on $(V,0)$, i.e. the factor of the module $\Omega^k_{\CC^N,0}$
of $k$-forms on $(\CC^N,0)$ by the submodule generated by
$f\cdot\Omega^k_{\CC^N,0}$ and $df\wedge\Omega^k_{\CC^N,0}$ for all $f$
from the ideal of functions vanishing on $(V,0)$.
Consider the Koszul complex $(\Omega^\bullet_{V,0}, X)$:
$$
0 \longleftarrow \cO_{V,0}
\stackrel{X}{\longleftarrow} \Omega^1_{V,0}
\stackrel{X}{\longleftarrow} ... \stackrel{X}{\longleftarrow}
\Omega^n_{V,0} \longleftarrow 0\,,
$$
where
the arrows are given by contraction with the vector field $X$.
The sheaves $\Omega^i_{V,0}$ are coherent sheaves and the homology groups
of the complex $(\Omega^\bullet_{V,0}, X)$ are concentrated at
the origin and therefore are finite dimensional. The following
definition is due to X.~G\'omez-Mont
\cite{GomezMont98}.

\begin{definition} The {\em homological index}
$\,{\ind}_{\rm hom}(X; V,0)={\ind}_{\rm hom}\, X$ of the vector field
$X$ on $(V, 0)$ is the Euler characteristic of the above
complex:
\begin{equation}\label{eq0}
{\ind}_{\rm hom}(X; V,0) ={\ind}_{\rm hom}\, X = \sum_{j=0}^n
(-1)^j h_j(\Omega^\bullet_{V,0},X)\,,
\end{equation}
where $h_j(\Omega^\bullet_{V,0},X)$ is the dimension of the
corresponding homology group as a vector space over $\bC$.
\end{definition}

The homological index satisfies the law of conservation of number
\cite[Theorem~1.2]{GomezMont98}. In the case when $V$ is a hypersurface,
G\'omez-Mont has shown that the homological index is equal to the
GSV index.

Given a holomorphic 1-form $\omega$ on $(V,0)$ with an isolated singularity,
we consider the complex $(\Omega^\bullet_{V,0}, \wedge\omega)$:
$$
0 \longrightarrow \cO_{V,0} \stackrel{\wedge \omega}{\longrightarrow}
\Omega^1_{V,0}
\stackrel{\wedge \omega}{\longrightarrow} ...
\stackrel{\wedge \omega}{\longrightarrow} \Omega^n_{V,0} \longrightarrow 0\,,
$$
where the arrows are given by the exterior product by the form $\omega$.

This complex is the dual of the Koszul complex considered above.
It was used by G.M.~Greuel in \cite{Greuel75} for complete intersections.
In \cite{EGS04} the definition of G\'omez-Mont was adapted to this case.

\begin{definition} The {\em homological index}
$\,{\ind}_{\rm hom}(\omega; V,0)={\ind}_{\rm hom}\, \omega$ of the 1-form
$\omega$ on $(V, 0)$ is $(-1)^n$ times the Euler characteristic of the above
complex:
\begin{equation}\label{eq1}
{\ind}_{\rm hom}(\omega; V,0) ={\ind}_{\rm hom}\, \omega= \sum_{j=0}^n
(-1)^{n-j} h_j(\Omega^\bullet_{V,0},\wedge\omega)\,,
\end{equation}
where $h_j(\Omega^\bullet_{V,0},\wedge\omega)$ is the dimension of the
corresponding homology group as a vector space over $\bC$.
\end{definition}

In \cite{EGS04}, there was proved the following statement.

\begin{theorem} \label{t:homological} Let $\omega$ be a holomorphic
1-form on $V$ with an isolated singularity at the origin $0$.

\medskip
\item {\rm(i)} If $V$ is smooth, then ${\ind}_{\rm hom}\, \omega$
equals the usual index of the holomorphic 1-form $\omega$.

\medskip
\item {\rm(ii)} The homological index satisfies the law of
conservation of number: if
$\omega'$ is a holomorphic 1-form on $V$ close to $\omega$, then:
$$
{\ind}_{\rm hom}(\omega; V,0) = {\ind}_{\rm hom}(\omega'; V,0)
+ \sum {\ind}_{\rm hom}(\omega';V,x) \,,
$$
where the sum on the right hand side is over all those points $x$
in a small punctured neighbourhood of the origin $0$ in $V$ where the
form $\omega'$ vanishes.

\medskip
\item {\rm(iii)} If $(V,0)$ is an isolated complete intersection
singularity, then the homological index ${\ind}_{\rm hom}\,
\omega$ coincides
with the GSV index ${\ind}_{\rm GSV}\,\omega$.
\end{theorem}

Statement (ii) follows from \cite{GG02}.

\begin{remark}
We notice that one has
an invariant for functions on $(V, 0)$ with an isolated singularity at the
origin defined by $f \mapsto {\ind}_{\rm hom}\,df$. By the theorem above, if
$(V,0)$ is an isolated complete intersection singularity, this invariant
counts the number of critical points of the function $f$ on a Milnor fibre
of the {\icis} $V$.
\end{remark}

\begin{remark}\label{remMvS}
Let $(C,0)$ be a curve singularity and let $(\bar{C},\bar{0})$
be its normalization. Let $\tau=\dim \rm {Ker}(\Omega^1_{C,0}\to
\Omega^1_{\bar{C},\bar{0}})$, $\lambda=\dim \omega_{C,0}/
c(\Omega^1_{C,0})$, where $\omega_{C,0}$ is the dualizing
module of Grothendieck, $c:\Omega^1_{C,0}\to\omega_{C,0}$
is the class map (see \cite{BG80}). In \cite{MvS01} there is considered
a Milnor number of a function $f$ on a curve singularity
introduced by V.Goryunov \cite{Goryunov00}. One
can see that this Milnor number can be defined for a 1-form
$\omega$ with an isolated singularity on $(C,0)$ as
well (as $\dim \omega_{C,0}/\omega\wedge \cO_{C,0}$) and is
equal to ${\ind}_{\rm hom}\,\omega + \lambda - \tau$.
\end{remark}

The laws of conservation of numbers for the homological and the radial
indices of 1-forms together with the fact that these two indices
coincide on smooth varieties imply that their difference is a
locally constant, and therefore constant, function on the space
of 1-forms on $V$ with isolated singularities at the origin.
Therefore one has the following statement.

\begin{proposition}\label{propMilnor} Let $(V,0)$ be a germ of a complex analytic
space of pure dimension $n$ with an isolated singular point at
the origin. Then the difference
$${\ind}_{\rm hom}\,\omega - {\ind}_{\rm rad}\,\omega$$
between the homological and the radial index
does not depend on the 1-form $\omega$.
\end{proposition}

If $(V,0)$ is an {\icis}, this difference is equal to the Milnor number
of $(V,0)$.
Together with this fact proposition~\ref{propMilnor} permits to
consider the difference
$$
\nu(V,0): = {\ind}_{\rm hom}(\omega;V,0) - {\ind}_{\rm rad}(\omega;V,0)
$$
as a generalized Milnor number of the singularity $(V,0)$.

There are other invariants of isolated singularities of complex analytic
varieties which coincide with the Milnor number for isolated complete
intersection singularities. One of them is $(-1)^n$ times the reduced
Euler characteristic (i.e., the Euler characteristic minus 1) of the
absolute de Rham complex of $(V, 0)$:
$$0 \longrightarrow {\mathcal O}_{V,0} \stackrel{d}{\longrightarrow}
\Omega^1_{V,0} \stackrel{d}{\longrightarrow} \ldots
\stackrel{d}{\longrightarrow} \Omega^n_{V,0} \longrightarrow 0\,:$$
$$
\bar{\chi}(V,0):=\sum_{i=0}^n (-1)^{n-i} h_i(\Omega^\bullet_{V,0},d)-(-1)^n.
$$

In \cite{EGS04} it was shown that:

\begin{theorem} \label{thm:EGS} One has
$$
\nu(V, 0)=\bar{\chi}(V,0)
$$
if
\begin{itemize}
\item[{\rm (i)}] $(V,0)$ is a curve singularity,
\item[{\rm (ii)}] $(V,0) \subset (\CC^{d+1},0)$ is the cone over the
rational normal curve in
$\CC\PP^d$.
\end{itemize}
\end{theorem}

Statement (i) implies that the invariant $\nu(V,0)$ is different from
the Milnor number introduced by R.-O.~Buchweitz and
G.-M.~Greuel \cite{BG80} for curve singularities. Statement (ii) was
obtained with the help of H.-Ch. Graf von
Bothmer and R.-O.~Buchweitz. For $d=4$ this is Pinkham's example
\cite{Pinkham74} of a singularity which has
smoothings with different Euler characteristics.


\section{Euler obstruction}\label{sec5}

The idea of the Euler obstruction emerged from \cite{MacPherson74}
where the Euler obstruction of a singular point of a complex
analytic variety was defined (the definition was formulated in
terms of obstruction theory in \cite{BS81}).

The Euler obstruction of an isolated singular point (zero) of a
vector field on a (singular) complex analytic variety was
essentially defined in \cite{BMPS04} (though formally speaking
there it was defined only for holomorphic functions, i.e. for
corresponding gradient vector fields). The main result of
\cite{BMPS04} (see Theorem~\ref{theo:BMPS} below) gives a relation
between the Euler obstruction of a function germ on a singular
variety $V$ and the Euler obstructions of singular points of $V$
itself.

In \cite{BMPS04}, there is introduced the notion of the local Euler
obstruction of a holomorphic function
with an isolated critical point on the germ of a complex analytic
variety.
It is defined as follows.

Let $(V,0) \subset (\CC^N,0)$ be the germ of a purely $n$-dimensional
complex analytic variety with a Whitney stratification
$V=\bigcup_{i=1}^q V_i$ and
let $f$ be a holomorphic function defined in a neighbourhood of the origin
in $\CC^N$ with an isolated singular point on $V$ at the origin. Let $\eps>0$
be small enough such that the function $f$ has no singular points on
$V\setminus\{0\}$ inside the ball $B_\eps$. Let ${\rm grad}\,f$ be
the gradient vector field of $f$ as defined in Section~\ref{sec3}.
Since $f$ has no singular
points on $V \setminus \{0\}$ inside the ball $B_\eps$, the angle of
${\rm grad}\, f(x)$ and the tangent space $T_xV_i$ to a point
$x \in V_i \setminus \{ 0 \}$ is less than $\pi/2$. Denote by
$\zeta_i(x) \neq 0$ the projection of ${\rm grad}\, f(x)$ to the tangent
space $T_xV_i$.

Following the construction in \cite{BMPS04}, the vector fields $\zeta_i$
can be glued together to obtain a stratified vector field ${\rm grad}_V f$
on $V$ such that ${\rm grad}_V f$ is homotopic to the restriction of
${\rm grad}\,f$ to $V$ and satisfies ${\rm grad}_V f(x) \neq 0$ unless $x=0$.

Let $\nu: \widehat{V} \to V$ be the Nash transformation of the variety $V$
defined as follows. Let $G(n,N)$ be the Grassmann manifold of
$n$-dimensional vector subspaces of $\CC^N$. For a suitable neighbourhood
$U$ of the origin in $\CC^N$, there is a natural map
$\sigma: V_{\rm reg} \cap U \to U \times G(n,N)$ which sends a point $x$
to $(x, T_x V_{\rm reg})$ ($V_{\rm reg}$ is the non-singular part of $V$).
The Nash transform $\widehat{V}$ is the closure of the image
${\rm Im}\,\sigma$ of the map $\sigma$ in $U \times G(n,N)$. The Nash
bundle $\widehat{T}$ over $\widehat{V}$ is a vector bundle of rank $n$
which is the pullback of the tautological bundle on the Grassmann manifold
$G(n,N)$. There is a natural lifting of the Nash transformation to a bundle
map from the Nash bundle $\widehat{T}$ to the restriction of the tangent
bundle $T\CC^N$ of $\CC^N$ to $V$. This is an isomorphism of $\widehat{T}$
and $TV_{\rm reg} \subset T\CC^N$ over the regular part $V_{\rm reg}$ of $V$.

The vector field ${\rm grad}_V f$ gives rise to a section $\widehat{\zeta}$
of the Nash bundle $\widehat{T}$ over the Nash transform $\widehat{V}$
without zeros outside of the preimage of the origin.

\begin{definition} {\rm \cite{BMPS04}} The {\em local Euler obstruction}
${\rm Eu}_{V,0} \, f$ of the function $f$ on $V$ at the origin is the
obstruction to extend the non-zero section $\widehat{\zeta}$ from the
preimage of a neighbourhood of the sphere $S_\eps= \partial B_\eps$ to
the preimage of its interior, more precisely its value $($as an element of
$H^{2n}(\nu^{-1}(V\cap B_\eps),\nu^{-1}(V\cap S_\eps))$\,$)$ on the fundamental
class of the pair $(\nu^{-1}(V\cap B_\eps), \nu^{-1}(V \cap S_\eps))$.
\end{definition}

The word {\em local} will usually be omitted.

\begin{remark} The local Euler obstruction can also be defined for the
restriction to $V$ of a real analytic function on $\RR^{2N}$. In particular,
one can take for $f$ the squared distance on $V$ to the origin. In this
case, the invariant ${\rm Eu}_{V,0} \, f$ is the usual local Euler
obstruction ${\rm Eu}_V(0)$ of the variety $(V,0)$
defined in \cite{MacPherson74, BS81, BLS00}.
\end{remark}

Denote by $M_f=M_{f,t_0}$ the Milnor fibre of $f$, i.e.\ the
intersection $V\cap B_\eps\cap f^{-1}(t_0)$ for a regular value $t_0$ of $f$
close to $0$. In \cite[Theorem~3.1]{BMPS04} the following result is proved.

\begin{theorem} \label{theo:BMPS}
Let $f: (V,0) \to (\CC,0)$ have an isolated singularity at $0 \in V$. Then
$$
{\rm Eu}_V(0) = \left( \sum_{i=1}^q
\chi(M_f \cap V_i) \cdot {\rm Eu}_V(V_i) \right) + {\rm Eu}_{V,0} \, f,
$$
where ${\rm Eu}_V(V_i)$ is the value of the Euler obstruction of $V$
at any point of $V_i$, $i=1, \ldots, q$.
\end{theorem}

One can also define the local Euler obstruction of a stratified vector field on a germ of a complex analytic variety $(V,0)$ with a Whitney stratification \cite{BSS05b}. Let $V_i$ be a stratum of the Whitney stratification with $0 \in V_i$ and $X_i$ be a vector field on $V_i$ with an isolated singularity at $0$. The Proportionality Theorem of \cite{BS81, BSS05b} states that the local Euler obstruction of a radial extension $X$ of $X_i$ at $0$ is equal to the local Euler obstruction ${\rm Eu}_V(0)$ of $(V,0)$ times the radial index ${\rm ind}_{\rm rad} \, (X;V,0)$ of $X$.

In \cite{EGGD} the definition of the local Euler obstruction
of a function was adapted to the case of a 1-form.
Let $\omega$ be a 1-form on a neighbourhood of the origin in $\CC^N$
with an isolated singular point on $V$ at the origin. Let $\eps>0$ be
small enough such that the 1-form $\omega$ has no singular points on
$V \setminus \{0\}$ inside the ball $B_\eps$.

The 1-form $\omega$ gives rise to a section $\widehat{\omega}$ of the
dual Nash bundle $\widehat{T}^\ast$ over the Nash transform $\widehat{V}$
without zeros outside of the preimage of the origin.

\begin{definition} The {\em local Euler obstruction} ${\rm Eu}_{V,0}\,\omega$
of the 1-form $\omega$ on $V$ at the origin is the obstruction to extend
the non-zero section $\widehat{\omega}$ from the preimage of a
neighbourhood of the sphere $S_\eps= \partial B_\eps$ to the preimage
of its interior, more precisely its value (as an element of the cohomology
group $H^{2n}(\nu^{-1}(V\cap B_\eps), \nu^{-1}(V \cap S_\eps))$\,) on the
fundamental class of the pair
$(\nu^{-1}(V\cap B_\eps), \nu^{-1}(V \cap S_\eps))$.
\end{definition}

\begin{remark} The local Euler obstruction can also be defined for a real
1-form on the germ of a real analytic variety if the last one is
orientable in an appropriate sense.
\end{remark}

\begin{example} Let $\omega = df$ for the germ $f$ of a holomorphic
function on $(\CC^N,0)$. Then ${\rm Eu}_{V,0} \, df$ differs from the
Euler obstruction ${\rm Eu}_{f,V}(0)$ of the function $f$ by the sign
$(-1)^n$. The reason is that for the germ of a holomorphic function
with an isolated critical point on $(\CC^n,0)$ one has
${\rm Eu}_{f,V}(0) = (-1)^n \mu_f$ (see \cite[Remark~3.4]{BMPS04}), whence
${\rm Eu}_{V,0} \, df= \mu_f$ ($\mu_f$ is the Milnor number of the germ $f$).
E.g., for $f(x_1, \ldots , x_n)=x_1^2 + \ldots + x_n^2$ the obstruction
${\rm Eu}_{f,V}(0)$ is the index of the vector field
$\sum_{i=1}^n \overline{x}_i \partial/ \partial x_i$ (which is equal
to $(-1)^n$), but the obstruction ${\rm Eu}_{V,0} \, df$ is the index of
the (holomorphic) 1-form $\sum_{i=1}^n x_i dx_i$ which is equal to 1.
\end{example}

The Euler obstruction of a 1-form satisfies the law of conservation of
number (just as the radial index). Moreover, on a smooth variety the
Euler obstruction and the radial index coincide. This implies the
following statement (cf.\ Theorem~\ref{theo:BMPS}). We set
$\overline{\chi}(Z):=\chi(Z)-1$ and call it the {\em reduced} (modulo
a point) Euler characteristic of the topological space $Z$ (though,
strictly speaking, this name is only correct for a non-empty space $Z$).

\begin{proposition} \label{prop4}
Let $(V,0) \subset (\CC^N,0)$ have an isolated singular point at the
origin and let $\ell: \CC^N \to \CC$ be a generic linear function. Then
$$
{\ind}_{\rm rad} \, (\omega;V,0) - {\rm Eu}_{V,0} \, \omega =
{\ind}_{\rm rad} \,(d\ell;V,0) = (-1)^{n-1}\overline{\chi}(M_\ell),
$$
where $M_\ell$ is the Milnor fibre of the linear function $\ell$ on $V$.
In particular $$
{\rm Eu}_{V,0} \, df = (-1)^n (\chi(M_\ell) - \chi(M_f)).
$$
\end{proposition}

An analogue of the Proportionality Theorem for 1-forms was proved in \cite{BSS05a}.

Now let $(V,0) \subset (\CC^N,0)$ be an arbitrary germ of an analytic
variety with a Whitney stratification $V=\bigcup_{i=0}^q V_i$,
where we suppose that $V_0=\{0\}$.
For a stratum $V_i$, $i=0, \ldots , q$, let $N_i$ be the normal slice in
the variety $V$ to the stratum $V_i$ ($\dim N_i = \dim V- \dim V_i$) at
a point of the stratum $V_i$ and let $n_i$ be the index of a generic
(non-vanishing) 1-form $d\ell$ on $N_i$:
$$
n_i = (-1)^{\dim N_i-1} \overline{\chi}(M_{\ell|_{N_i}}).
$$
In particular for an open stratum $V_i$ of $V$, $N_i$ is a point and
$n_i=1$. The strata $V_i$ of $V$ are partially ordered: $V_i \prec V_j$
(we shall write $i \prec j$) iff $V_i \subset \overline{V_j}$ and
$V_i \neq V_j$; $i \preceq j$ iff $i \prec j$ or $i=j$. For two strata
$V_i$ and $V_j$ with $i \preceq j$, let $N_{ij}$ be the normal slice of
the variety $\overline{V_j}$ to the stratum $V_i$ at a point of it
($\dim N_{ij}= \dim V_j - \dim V_i$, $N_{ii}$ is a point) and let $n_{ij}$
be the index of a generic
1-form $d\ell$ on $N_{ij}$:
$n_{ij}=(-1)^{\dim N_{ij}-1}\, \overline{\chi}(M_{\ell\vert_{N_{ij}}})$,
$n_{ii}=1$. Let us define the Euler obstruction ${\rm Eu}_{Y,0} \, \omega$
to be equal to 1 for a zero-dimensional variety $Y$ (in particular
${\rm Eu}_{\overline{V_0},0}\,\omega =1$, ${\rm Eu}_{N_{ii},0}\,\omega =1$).
In \cite{EGGD}, the following statement was proved.

\begin{theorem} \label{theo4}
One has
$$
{\ind}_{\rm rad} \, (\omega;V,0) =
\sum_{i=0}^q n_i \cdot {\rm Eu}_{\overline{V_i},0}\, \omega.
$$
\end{theorem}

To write an "inverse" of the formula of Theorem~\ref{theo4}, suppose that
the variety $V$ is irreducible and $V=\overline{V_q}$. (Otherwise one can
permit $V$ to be reducible, but also permit the open stratum $V_q$ to be not
connected and dense; this does not change anything in Theorem~\ref{theo4}.)
Let $m_{ij}$ be the (M\"obius) inverse of the function $n_{ij}$ on the
partially ordered set of strata, i.e.
$$
\sum_{i \preceq j \preceq k} n_{ij}m_{jk} = \delta_{ik}.
$$
For $i \prec j$ one has
\begin{eqnarray*}
m_{ij} & = & \sum_{i=k_0 \prec k_1 \prec \ldots \prec k_r =j} (-1)^{r}
n_{k_0k_1}n_{k_1k_2} \ldots n_{k_{r-1}k_r}\\
& = &
(-1)^{\dim V -\dim V_i} \sum_{i=k_0 \prec \ldots \prec k_r =j}
\overline{\chi}(M_{\ell\vert_{N_{k_0k_1}}}) \cdot \ldots \cdot
\overline{\chi}(M_{\ell\vert_{N_{k_{r-1}k_r}}})\ .
\end{eqnarray*}

\begin{corollary}
One has
$$
{\rm Eu}_{V,0}\, \omega =
\sum_{i=0}^q m_{iq} \cdot {\ind}_{\rm rad} \, (\omega;\overline{V_i},0).
$$
In particular
\begin{eqnarray*}
\lefteqn{{\rm Eu}_{V,0}\, df = (-1)^{\dim V -1} \times}\\
&& \left( \overline{\chi}(M_{f\vert_V}) + \sum_{i=0}^{q-1}
\overline{\chi}(M_{f\vert_{\overline{V_i}}})
\sum_{i=k_0 \prec \ldots \prec k_r =q}
\overline{\chi}(M_{\ell\vert_{N_{k_0k_1}}}) \ldots
\overline{\chi}(M_{\ell\vert_{N_{k_{r-1}k_r}}}) \right).
\end{eqnarray*}
\end{corollary}

In \cite{LT81} the local Euler obstruction of a germ of a variety is related to polar invariants of the germ. A global version of the Euler obstruction of a variety and its
connections with global polar invariants are discussed in \cite{STV05a}.

\section{Algebraic, analytic, and topological formulae for indices}
\label{sec6} 
In the introduction, there were mentioned algebraic formulae for the
index of an analytic vector field or a 1-form on a smooth manifold:
as the dimension of a ring in the complex setting and as the
signature of a quadratic form in the real one. It is natural to try
to look for analogues of such formulae for
vector fields or 1-forms on singular varieties. This appeared to be a
rather complicated problem. Other sorts
of formulae: analytic (usually as certain residues) or topological
ones are of interest as well.

A substantial progress in this direction was achieved for the GSV index.

Let $(V,0)$ be the germ of a hypersurface in $(\CC^{n+1}, 0)$
defined by a germ of a holomorphic function
$f: (\CC^{n+1},0) \to (\CC,0)$ with an isolated singular point at $0$. Let
$$X = \sum_{i=1}^{n+1} X_i \frac{\partial}{\partial x_i}$$
be a holomorphic vector field on $\CC^{n+1}$ tangent to $V$ with an
isolated (in $\CC^{n+1},0$) zero at the origin. This implies
that $Xf= h f$ for some $h \in {\mathcal O}_{\CC^{n+1},0}$. Define the
following ideals in ${\mathcal O}_{\CC^{n+1},0}$:
\begin{eqnarray*}
J_f &:= & (f_1, \ldots, f_{n+1}),\\
J_X & := & (X_1, \ldots, X_{n+1}), \\
J_1 & := & (h, J_X), \\
J_2 & := & (f, J_f), \\
J_3 & := & (f, J_X).
\end{eqnarray*}
X.~G\'omez-Mont \cite{GomezMont98} has proved the following formula
for the GSV index:

\begin{theorem}
For a vector field $X$ with an isolated zero in the ambient space tangent
to the germ of a hypersurface $(V,0)$ with an
isolated singularity at 0, we have: for even $n$
$${\rm Ind}_{\rm GSV} \, (X; V,0) = \dim {\mathcal O}_{\CC^{n+1},0}/J_X +
\dim {\mathcal O}_{\CC^{n+1},0}/J_1 + \dim
{\mathcal O}_{\CC^{n+1},0}/J_2\,,$$
for $n$ odd
$${\rm Ind}_{\rm GSV} \, (X; V,0) = \dim {\mathcal O}_{\CC^{n+1},0}/J_2 +
\dim {\mathcal O}_{\CC^{n+1},0}/J_3\,.$$
\end{theorem}

In this case the GSV index coincides with the homological index.
O.~Klehn \cite{Klehn03} generalized this formula for the homological
index to the case when the vector field has an isolated zero on the
hypersurface singularity $(V,0)$ but not necessarily in $(\CC^N,0)$.
Recently, H.-Ch.~Graf von Bothmer, X.~G\'omez-Mont and the first author
gave formulae to compute the homological index in the case when $V$
is a complete intersection \cite{BEG05}.

In \cite{GM97,GM99} algebraic formulae for the index of a real
analytic vector field with an algebraically
isolated singular point at the origin tangent to a real hypersurface
with an algebraically isolated singularity
at 0 are derived. The index is expressed as the signature of a
certain non-degenerate quadratic form for an
even-dimensional hypersurface and as the difference of the signatures
of two non-degenerate quadratic forms in
the odd-dimensional case.

O.~Klehn \cite{Klehn05} proved that the GSV index coincides with the
dimension of a certain explicitly constructed vector space in the case
when $(V,0)$ is an isolated complete intersection singularity of dimension
1 and $X$ is a vector field tangent to $V$ with an isolated zero on $V$
which is deformable in a certain way. He also obtained a signature formula
for the real GSV index in the corresponding real analytic case
generalizing the Eisenbud-Levine-Khimshiashvili formula.

Let $(V,0)\subset (\CC^{n+k},0)$ be an isolated complete intersection
singularity defined by an analytic map
$$
f = (f_1,...,f_k) : (\CC^{n+k},0) \to (\CC^k,0) \,.
$$
Let $X$ be a holomorphic vector field on $\CC^{n+k}$ defined in a
suitable neighbourhood $U$ of the origin, tangent to $V$, with an isolated
singular point at the origin. Let $C$ be the $k \times k$ matrix whose
entries are holomorphic functions on $U$ such that $X \cdot f = C f$.
Assume that $(x_1, \ldots , x_{n+k})$ is a system of coordinates on $U$
such that when $X$ is written as
$$
X = \sum_{i=1}^{n+k} X_i \frac{\partial}{\partial x_i},
$$
the sequence $(X_1, \ldots , X_n, f_1, \ldots, f_k)$ is regular. Let
$J$ denote the Jacobian matrix
$$
J= \left( \frac{\partial X_i}{\partial x_j} \right) .
$$
We denote by $c_q$ the coefficient
at $t^q$ in the formal power series expansion of
$$
\det \left( I_{n+k} - t \frac{\sqrt{-1}}{2 \pi} J \right) \left[
\det \left( I_k - t \frac{\sqrt{-1}}{2 \pi} C \right) \right]^{-1}
$$
in $t$, where $I_{n+k}$ and $I_k$ denote the identity matrices of
sizes $n+k$ and $k$ respectively. In \cite{LSS95} the following
formula for the GSV index is proved.

\begin{theorem} \label{thmLSS}
Let $\eps$ be small enough such that the real hypersurfaces
$|X_i|=\eps$ are in general position and let $Z$ be the set
$\{ f=0, |X_i| = \eps, 1 \leq i \leq n \}$. We assume that $Z$
is oriented so that
$$
d(\arg X_1) \wedge d(\arg X_2) \wedge \ldots \wedge d(\arg X_n)
$$
is positive. Then
$$
{\rm Ind}_{\rm GSV} \, (X; V,0) = \int_Z \frac{c_n dx_1 \wedge
\ldots \wedge dx_n}{\prod_{i=1}^n X_i}.
$$
\end{theorem}

For generalizations of Theorem~\ref{thmLSS} and related results see \cite{LS95, Suwa95, Suwa98, Suwa02, IS03}.

Now let $\omega$ be a holomorphic 1-form on the {\icis}
$(V,0)=f^{-1}(0)\subset(\CC^{n+k},0)$,
$f=(f_1, \ldots, f_k): (\CC^{n+k},0)\to(\CC^k,0)$, i.e. the
restriction to $(V, 0)$ of a holomorphic 1-form
$$\omega = \sum_{i=1}^{n+k}A_i(x)dx_i$$
on $(\CC^{n+k}, 0)$. Assume that $\omega$ has an isolated singular
point at the origin (on $(V,0)$). Let $I$ is the ideal
generated by $f_1$, ..., $f_k$ and the $(k+1) \times (k+1)$-minors of
the matrix
$$\begin{pmatrix} \frac{\partial f_1}{\partial x_1} & {\cdots } &
   \frac{\partial f_1}{\partial x_{n+k}} \\ {\vdots} & {\cdots} &
{\vdots} \\ \frac{\partial
f_k}{\partial x_1} & \cdots &
   \frac{\partial f_k}{\partial x_{n+k}} \\
A_1 & \cdots & A_{n+k}
\end{pmatrix}\,.
$$
Then one has the following formula for the GSV index \cite{EG01, EG03a, EG03b}.

\begin{theorem}\label{sec6theoAlg}
One has
$$
{\rm ind}_{\rm GSV} \, (\omega;V,0) = \dim {\mathcal O}_{\CC^{n+k},0}/I.
$$
\end{theorem}

(Note that there is a minor mistake in the proof of this theorem in
\cite{EG03a} which is corrected in
\cite{EG05b}.) This formula was obtained by L\^e D.T. and G.-M.
Greuel for the case when $\omega$ is the
differential of a function (\cite {Greuel75, Le74}).

T.~Gaffney \cite{Gaffney05} described connections between the GSV index of a holomorphic 1-form on an {\icis} and the multiplicity of pairs of certain modules.

In \cite{EGMZ}, there was constructed a quadratic form on the algebra
${\mathcal A}:= {\mathcal O}_{\CC^{n+k},0} /I$
generalizing the Eisenbud-Levine-Khimshiashvili quadratic form
defined for a smooth $V$. This is defined as
follows.

Let $F:(\CC^{n+k}\times\CC_\eps^M,0)\to
(\CC^k\times\CC_\eps^M,0)$ be a deformation of the map
$f:(\CC^{n+k},0)\to(\CC^k,0)$
($F(x, \eps)=(f_\eps(x), \eps)$, $f_0=f$, $f_\eps=(f_{1\eps}, \ldots,
f_{k\eps})$) and let $\omega_\eps$ be a deformation of the form $\omega$
(defined in a neighbourhood of the origin in $\CC^{n+k}\times\CC_\eps^M$) such
that, for generic
$\eps\in (\CC_\eps^M,0)$, the preimage $f_\eps^{-1}(0)$ is smooth and
(the restriction of) the form $\omega_\eps$ to it has only non degenerate
singular
points.

Let $\Sigma\subset (\CC^M_\eps,0)$ be the germ of the set of the values of
the parameters $\eps$ from $(\CC^M_\eps,0)$ such that either the preimage
$f_\eps^{-1}(0)$ is singular or the restriction of the 1-form $\omega_\eps$
to it has degenerate singular points.

Let $\sum\limits_{i=1}^n  A_i dy^i$ be a 1-form on a smooth complex
analytic manifold of dimension $n$ (in local coordinates $y^1$, \dots,
$y^n$).
At
points $P$ where the 1-form $\omega$ vanishes, the Jacobian matrix
$\mathcalJ=\left(\frac{\partial A_i}{\partial y^j}\right)$
defines the tensor
$\sum\limits_{i,j}\frac{\partial A_i}{\partial y^j}\, dy^i\otimes dy^j$
of type $(0, 2)$ (i.e. a bilinear form on the tangent space). The determinant
$J$ of the matrix $\mathcalJ=\left(\frac{\partial A_i}{\partial
y^j}\right)$ is not
a scalar (it depends on the choice of local coordinates). Under a change of
coordinates it is multiplied by the square of the Jacobian of the change
of coordinates (since the Jacobian matrix $\mathcalJ$ is transformed to
$C^T\mathcalJ C$).
Therefore it should be considered as the coefficient in the tensor
$$
J(dy^1\wedge\ldots\wedge dy^m)^{\otimes 2}
$$
of type $(0,2m)$. In this sense the Jacobian of a 1-form is a sort of a
"quadratic differential". To get a number $\widetilde{J}(P)$, one can
divide this tensor by the
tensor square of a volume form.

Let us fix volume forms on $\CC^{n+k}$ and $\CC^k$, say, the
standard ones
$\sigma_{n+k}=dx_1\wedge\ldots\wedge dx_{n+k}$ and
$\sigma_k=dz_1\wedge\ldots\wedge
dz_k$
where $x_1$, \dots, $x_{n+k}$ and $z_1$, \dots, $z_k$ are Cartesian coordinates
in $\CC^{n+k}$
and in $\CC^k$ respectively. There exists (at least locally) an
$n$-form $\sigma$
on $\CC^n$ such that $\sigma_{n+k}=f^*\sigma_k\wedge\sigma$. Let
$\eps \not\in \Sigma$. The restriction
of the form $\sigma$ to the manifold $f_\eps^{-1}(0)$ is well defined
and is
a volume form on $f_\eps^{-1}(0)$. We also denote it by $\sigma$.
Let $P_1$, \dots, $P_\nu$ be the (non-degenerate) singular points of the 1-form
$\omega_\eps$ on the $n$-dimensional manifold $f_\eps^{-1}(0)$.
For $\eps\notin\Sigma$ and a germ $\varphi\in\mathcalO_{\CC^n,0}$, let
\begin{equation}\label{Formula1}
R(\varphi,\eps):=\sum\limits_{i=1}^\nu \frac{\varphi(P_i)}{\widetilde
J_\eps(P_i)}.
\end{equation}
For a fixed $\varphi$ the function $R_\varphi(\eps):=R(\varphi,\eps)$
is holomorphic in
the complement of the bifurcation diagram $\Sigma$.

In \cite{EGMZ}, it is proved that the function $R_\varphi(\eps)$ has
removable singularities on the bifurcation
diagram $\Sigma$. This means that $R_\varphi(\eps)$ can be extended
to a holomorphic
function on $(\CC^M_\eps,0)$ (which we denote by the same symbol). Let
\begin{equation}
R(\varphi):=R_\varphi(0).
\end{equation}
This defines a linear function $R$ on ${\mathcal O}_{\CC^n,0}$. One
can show that this linear function vanishes on the
ideal $I$ and hence defines a linear function $R$ on the algebra
${\mathcal A}$. Using this function, we can define
a quadratic form
$$
Q(\varphi, \psi)=R(\varphi \cdot \psi)
$$
on the algebra ${\mathcal A}$ generalizing the
Eisenbud-Levine-Khimshiashvili quadratic form defined
for a smooth $V$.

A similar form can be defined on a module of differential forms
associated to the pair $(V,\omega)$, see \cite{EGMZ}.
A special case of such a quadratic form
(for curve singularities)
was considered in \cite{MvS90}.
In the smooth case this quadratic form can be identified with the
one on the algebra ${\mathcal A}$ by an isomorphism of the underlying
spaces. This is no longer true for an {\icis}. In general,
these two quadratic forms have even different ranks.

Let $g : (\CC^{n+1},0) \to (\CC,0)$ be a real analytic function (i.e.
a function taking real values on $\RR^{n+1} \subset \CC^{n+1}$) with
an isolated singular point at 0. Then there is the following topological
formula for the index of the gradient vector field ${\rm grad} \, g$ at
the origin. Let $\eta$ be a sufficiently small positive real number.
The complex conjugation induces involutions on the Milnor fibres
$M_{g,\eta}$ and $M_{g,-\eta}$ of $g$. Their actions on the homology groups
$H_n(M_{g,\eta})$ and $H_n(M_{g,-\eta})$ will be denoted by $\sigma_+$ and
$\sigma_-$ respectively. On these homology groups we have the intersection
forms $\langle \cdot , \cdot \rangle$. We consider the symmetric bilinear
forms $Q_+(a,b):=\langle \sigma_+ a,b \rangle$ and
$Q_-(a,b):= \langle \sigma_- a,b \rangle$ on $H_n(M_{g,\eta})$ and
$H_n(M_{g,-\eta})$ respectively. The following statement was conjectured
by Arnold \cite{Arnold78} and proved in \cite{GuseinZade84, Varchenko85}.

\begin{theorem}
If $n$ is even then
$$
{\rm ind}_0 \, {\rm grad} \,g = (-1)^{\frac{n}{2}} \frac{1}{2} (
{\rm sgn}\, Q_- - {\rm sgn} Q_+)\,,
$$
where ${\rm grad} \,g$ is the gradient vector field of the function $g$
on $\RR^{n+1}$.
\end{theorem}

If ${\rm Var} : H_n(M_{g,\eta}, \partial M_{g,\eta}) \to
H_n(M_{g,\eta})$ denotes the variation operator of the
singularity of $g$, then there is another formula for the index in
\cite{GuseinZade84}:

\begin{theorem}
One has
$${\rm ind}_0 \, {\rm grad} \,g = (-1)^{\frac{n(n-1)}{2}} {\rm sgn}
\, {\rm Var}^{-1} \sigma_+.$$
\end{theorem}

In \cite{EG99} a generalization of such a formula for the radial
index of a gradient vector field on an
algebraically isolated real analytic isolated complete intersection
singularity in $\CC^{n+k}$ was obtained.

In \cite{Seade95, SS96} formulas are given evaluating the GSV
index of a singular point of a vector field on an isolated
hypersurface or complete intersection singularity on a resolution
of the singularity.

A.~Esterov has extended the method of A.Khovanski of computation of
intersection numbers from complete intersections to determinantal
varieties. Using this he gave formulae for the GSV index of a 1-form
on an {\icis} and for its generalization to collections of 1-forms
(see Section~\ref{sec8}) in terms of Newton diagrams of the equations
of the {\icis} and of the components of the 1-forms. See \cite{Esterov05b}
for a short account of the results (a full version has been submitted
to a journal; see also \cite{Esterov05} for some preliminary results).

In \cite{Klehn02} a residue formula for the GSV index of a
holomorphic 1-form on an isolated surface
singularity in the spirit of Theorem~\ref{thmLSS} is proved.


\section{Indices of meromorphic 1-forms}\label{sec7} 

Nonzero holomorphic vector fields (as well as holomorphic 1-forms)
on a compact complex manifold rarely exist. Therefore it is interesting
to consider meromorphic vector fields or 1-forms of which there are a
lot (at least for projective manifolds). For a meromorphic vector field
defined by a holomorphic section of the vector bundle $TM \otimes L$ ($L$
is a holomorphic line bundle on $M$), the sum of the indices of its zeros
is equal to the corresponding characteristic number
$\langle c_n(TM\otimes L),[M]\rangle$
of the vector bundle $TM \otimes L$ (see, e.g., \cite{BB70}).

Let $\alpha$ be a meromorphic $1$-form on a compact complex manifold $M^n$
which means that $\alpha$ is a holomorphic $1$-form outside of a positive
divisor $D$ and in a neighbourhood of each point of $M$ the form $\alpha$
can be written as $\widehat\alpha/F$ where $F=0$ is a local equation of
the divisor $D$ and $\widehat\alpha$ is a holomorphic $1$-form. Let $L$
be the line bundle associated to the divisor $D$, i.e., $L$ has a holomorphic
section $s$ with zeros on $D$. Then $\omega=s\alpha$ is a holomorphic
section of the vector bundle $T^{\ast}M{\otimes}L$. In some constructions
(say, as in \cite{BB70}) one defines a meromorphic $1$-form on $M$ simply as
a holomorphic section of the tensor product $T^{\ast}M\otimes L$ for a
holomorphic line bundle $L$. This definition is somewhat different from the
one formulated above. E.g., in this case only the class of the divisor of
poles of a meromorphic $1$-form is defined, not the divisor itself.
Moreover, in this setting, the value of a meromorphic 1-form on a vector
field is not a function, but a section of the line bundle $L$. In the
sequel we use the notation $\omega$ for a section of the vector bundle
$T^{\ast}M\otimes L$, for short calling it a meromorphic 1-form as well.
Suppose that the section $\omega$ has isolated zeros. Then the sum of
their indices is equal to the characteristic number
$\langle c_n(T^{\ast}M\otimes L),[M]\rangle$ and thus depends on $L$.

For a meromorphic $1$-form on a smooth compact complex curve $M$ the
characteristic number $\langle c_1(T^{\ast}M),[M]\rangle=-\chi(M)$
is equal to the
number of zeros minus the number of poles counted with multiplicities.
Therefore to express the Euler characteristic of a compact manifold
in terms of singularities of a meromorphic $1$-form one has to take
the divisor $D$ of its poles into account as well. If a meromorphic
$1$-form on a manifold $M^n$ is defined simply as a section of
$T^{\ast}M\otimes L$, the divisor of poles is not defined and thus
one also cannot define singular points of the 1-form on its pole locus.

One would like to give a Poincar\'{e}-Hopf type formula for meromorphic
1-forms, i.e., to express the Euler characteristic of a compact manifold
or of a smoothing of a complete intersection in terms of singularities
of a meromorphic $1$-form. For that it is possible to introduce a suitable
notion of an index of a germ of a meromorphic $1$-form (with an additional
structure) on an {\icis}, so that the indices of the singular points sum
up to (plus-minus) the Euler characteristic of a smoothing. There are
several descriptions of the index of a meromorphic 1-form, one of
them being the alternating sum of the dimensions of certain algebras.

A Poincar\'e-Hopf type formula can be considered as one describing a
localization of an invariant (say, of the Euler characteristic) of a
manifold at singular points of, say, a vector field or a $1$-form,
i.e., a representation of the invariant as the sum of integer invariants
(''indices'') corresponding to singular points. Therefore one first has
to define singular points. Let $M^n$ be a complex manifold, let $L$ be
a line bundle on $M$, and let $\omega$ be a holomorphic section of
the bundle $T^{\ast}M\otimes L$. One has zeros of $\omega$ on $M$, but
the pole locus of $\omega$ is not well defined. In order to discuss
singular points of $\omega$ on its pole locus (which is necessary as we
saw in the example when $M$ was a curve), we have to fix this locus.
This means that we have to choose a holomorphic section $s=s_1$ of the
line bundle $L$ or to choose its zero divisor $D=D_1$. One can say that
we have to consider $\alpha=\omega/s$, i.e., to proceed with our initial
definition of a meromorphic $1$-form.

For the further setting we suppose that the divisor $D$ of poles of the
$1$-form $\omega$ is non-singular (in particular, reduced). This is not very
essential for this section, however, this makes the discussion simpler.
Since $D$ is a submanifold of $M$, there is a well defined map
$T^{\ast}M\raisebox{-0.5ex}{$\vert$}{}_{D}\rightarrow T^{\ast}D$ and thus
$(T^{\ast}M\otimes L)\raisebox{-0.5ex}{$\vert$}{}_{D}\rightarrow
T^{\ast} D\otimes L\raisebox{-0.5ex}{$\vert$}{}_{D}$ (the restriction
of meromorphic $1$-forms to $D$). (Using the 1-form $\alpha$, one can
also define only a section of the vector bundle
$T^{\ast} D\otimes L\raisebox{-0.5ex}{$\vert$}{}_{D}$, but not
a meromorphic 1-form on $D$ with precisely defined pole locus.)

Let $\omega_1$ be the restriction of $\omega$ to $D_1$. It is a holomorphic
section of the vector bundle
$T^{\ast}D\otimes L\raisebox{-0.5ex}{$\vert$}{}_{D}$. Its zeros are well
defined and should be considered as singular points of the meromorphic
$1$-form $\omega$ on the pole locus. To discuss its singular points on
its pole locus we again have to fix a divisor. Suppose that there exists
a (positive) divisor $D_2$ on $M$ which is the zero locus of another
section $s_2$ of the same line bundle $L$ and which intersects $D_1$
transversally (in particular, this means that $D_2$ is non-singular at its
intersection points with $D_1$ and $D_1\cap D_2$ is non-singular as well).
One has the restriction of $\omega$ to $D_1\cap D_2$ and its zeros there.

Going on
this way we arrive at the situation when we have fixed $n$
divisors $D_1,\ldots, D_n$ (zeros of sections $s_1$, \dots, $s_n$ of the
line bundle $L$) so that, for each $i=1,\ldots n$, $D_1\cap\ldots \cap D_i$
is non-singular. The set of singular points of the 1-form $\omega$ is the
union of the zeros of $\omega$ itself and of the restrictions of $\omega$
to $D_1\cap\ldots\cap D_i$ for all $i=1,\ldots, n$. (For $i=n$ the
intersection $D_1\cap\ldots\cap D_n$ is zero-dimensional and all its points
should be considered as zeros of the 1-form $\omega$.) One can say that
we have to consider a collection of meromorphic 1-forms $\omega/s_1$,
\dots , $\omega/s_n$ on $M$ proportional to each other.

There is the following Poincar\'e-Hopf type formula for meromorphic
$1$-forms. Let $M^n$ be a compact complex manifold, let $\omega$ be
a meromorphic $1$-form on $M$, that is, a holomorphic section of the bundle
$T^{\ast}M\otimes L$ where $L$ is a holomorphic line bundle with nonzero
holomorphic sections. Suppose that $D_1=D, D_2,\ldots, D_n$ are zero
divisors of holomorphic sections of the line bundle $L$ such that, for
each $i=1,\ldots,n$, $D_1\cap\ldots\cap D_i$ is non-singular. Suppose
that the form $\omega$ itself and its restrictions to the submanifolds
$D_1 \cap \ldots \cap D_i$, $i=1, \ldots , n$, have only isolated zeros.
Let $m_0$ (respectively $m_i$, $i=1, \ldots , n$) be the number of zeros
of the form $\omega$ (respectively, of the restriction of $\omega$ to
the intersection $D_1\cap\ldots\cap D_i$) counted with multiplicities.
In particular $m_n$ is the number of points in $D_1 \cap\ldots\cap D_n$.

\begin{theorem}\label{sec7theo1}
$c_n(T^{\ast}M)[M]=(-1)^n\chi(M^n)=m_0-m_1+\ldots +(-1)^n m_n.$
\end{theorem}

\begin{remark}
One can show that
$$
m_0=(-1)^n(\chi(M)-\chi(D))=(-1)^n\chi(M\setminus D).
$$
Suppose that all zeros of $\omega$ are outside of $D$ (one can consider
this situation as the generic one). Let $s$ be the holomorphic section
of the line bundle $L$ with zeros on $D$. Then $\alpha=\omega/s$ is a
holomorphic $1$-form on $M\setminus D$ with simple poles along $D$.
So in this case the number $m_0$ of zeros of the holomorphic $1$-form
$\alpha$ on $M\setminus D$ coincides with $(-1)^n\chi(M\setminus D)$.
This is the relation which holds for holomorphic $1$-forms on compact
manifolds ($M\setminus D$ is not compact).
\end{remark}

\begin{example}
Let
$$
\alpha=\frac{xdy-ydx+dz}{x^2+4y^2+z^2+1}
$$
be a meromorphic 1-form on the projective space $\CC\PP^3$ ($x$, $y$, $z$
are affine coordinates). One can see that the zeros of the corresponding
$\omega$ on $\CC\PP^3$ and also the zeros of $\omega|_D$, $D=D_1=
\{x^2+4y^2+z^2+1=0\}$, are isolated and $m_0=0$, $m_1=4$. To define
other singular points one has to choose $D_2$ and $D_3$ (e.g.,
$D_2=\{x^2+y^2+4z^2=0\}$, $D_3=\{x^2+y^2+z^2=0\}$). One has
$m_2=8$, $m_3=8$, $0-4+8-8=-4=(-1)^3\chi(\CC\PP^3)$. Note
that as a meromorphic 1-form on $\CC\PP^3$ with poles on the
hypersurface $P_d(x, y, z)=0$, $\deg P_d=d$, it is natural to take
$$
\frac{A_{d-2}dx+B_{d-2}dy+C_{d-2}dz}{P_d},
$$
where $A_{d-2}$, $B_{d-2}$, and $C_{d-2}$ are polynomials in $x$, $y$, $z$
of degree $d-2$ ($\alpha$ is not of this form). However, one can show
that such a 1-form has non-isolated zeros on $\CC\PP^3$ (at infinity).
\end{example}

Now let $V^n$ be a compact subvariety of a complex manifold $M^{n+k}$
such that in a neighbourhood of each point $V$ is a complete intersection
in $M$ with only isolated singularities ($M$ is not supposed to be
compact). Let $\omega$ be a meromorphic 1-form on $M$, i.e., a
holomorphic section of the bundle $T^{\ast}M\otimes L$, where $L$
is a line bundle on $M$ with nonzero holomorphic sections. By a smoothing
$\widetilde V$ of $V$ we understand a smooth ($C^\infty$) manifold
which is obtained from $V$ by smoothing its singular points (which are
{\icis}) in the usual way (as complex {\icis}). (If $V$ is a complete
intersection in the projective space $\CC\PP^{n+k}$, a smoothing of $V$
can be obtained as an analytic submanifold of $\CC\PP^{n+k}$ as well.)

It is possible to formulate a Poincar\'e-Hopf type formula for the Euler
characteristic of $V$ and of its smoothing $\widetilde V$ in terms
of singular points of the 1-form $\omega$ on $V$. As above we suppose
that there exist $n$ holomorphic sections $s_1$, \dots, $s_n$ of the
line bundle $L$ with zero divisors $D_1$, \dots, $D_n$ such that, for each
$i=1, \ldots, n$, the intersection $V\cap D_1\cap\ldots\cap D_i$ is of
dimension $n-i$ and has only isolated singularities (and thus is reduced
for $i<n$). We also suppose that the restriction of $\omega$ to the
non-singular part of $V$ and its restriction to the non-singular part
of $V\cap D_1\cap\ldots\cap D_i$ ($i=1, \ldots, n$) have only isolated
zeros. As singular points of the 1-form $\omega$ on $V$ we consider
all zeros of its restrictions to the non-singular parts of $V$ and of
the intersections $V\cap D_1\cap\ldots\cap D_i$ ($i=1, \ldots, n$)
and all singular points of $V$ and of $V\cap D_1\cap\ldots\cap D_i$
as well.

Let $P$ be a singular point of the 1-form $\omega$ on $V$. Suppose that the
divisors $D_1$, \dots, $D_\ell$ ($0\le \ell\le n$) pass through the point $P$
and (if $\ell<n$) $D_{\ell+1}$ does not. In a neighbourhood of the point $P$,
in some local coordinates on $M$ centred at the point $P$, we have the
following situation. The variety ({\icis}) $V$ (possibly non-singular) is
defined by $k$ equations $f_1=\ldots=f_k=0$. Let divisors $D_1$, \dots,
$D_\ell$ be defined by equations $f_{k+1}=0$, \dots, $f_{k+\ell}=0$. Let us
remind that $\{f_1=\ldots=f_k=f_{k+1}=\ldots=f_{k+i}=0\}$ is an {\icis} for
each $i=0, 1, \ldots, \ell$. After choosing a local trivialization of the
line bundle $L$, the 1-form $\omega$ can be written as
$\sum_{i=1}^{n+k}A_i(x)dx_i$ where $A_i(x)$ are holomorphic. One can say
that we consider a collection of meromorphic 1-forms $\omega/f_{k+1}$,
\dots, $\omega/f_{k+\ell}$ on $V$ proportional to each other. Let us denote
the set of local data $(V, \omega, f_{k+1}, \ldots, f_{k+\ell})$ by $\Omega$.

Let ${\ind}^{(0)}\Omega$ (respectively ${\ind}^{(i)}\Omega$,
$i=1, \ldots, \ell$) be the index of the holomorphic 1-form $\omega$ on
$V$ (respectively on the {\icis} $V\cap D_1\cap\ldots\cap D_i$) defined
in Section~\ref{sec3}.

\begin{definition}
The alternating sum
$$
{\ind}_P\Omega={\ind}^{(0)}\Omega-
{\ind}^{(1)}\Omega+\ldots+(-1)^\ell {\ind}^{(\ell)}\Omega
$$
is called the index of the 1-form $\omega$ at the point $P$ with respect
to the divisors $D_1$, \dots, $D_n$.
\end{definition}

According to the statements from Section~\ref{sec3} (see also
Section~\ref{sec6} for the last point) one has the following
three equivalent descriptions of the index ${\ind}_P\Omega$.

{\bf 1.} Let $B_\delta$ be the ball of sufficiently small
radius $\delta$ centred at the origin in $\CC^{n+k}$ and let $S_\delta$
be its boundary. Let $\eps=(\eps_1, \ldots, \eps_k, \ldots, \eps_{k+\ell})
\in(\CC^{k+\ell}, 0)$ be small enough and such that $\eps^{(k+i)}=
(\eps_1, \ldots, \eps_k, \ldots, \eps_{k+i})$ is not a critical value
of the map $F_{k+i}=(f_1, \ldots, f_k, \ldots, f_{k+i}):(\CC^{n+k},0)
\to(\CC^{k+i},0)$ for each $i=0, 1, \ldots, \ell$. The restriction of
the 1-form $\omega$ to the smooth manifold $V_i=F_{k+i}^{-1}(\eps^{(k+i)})$
has isolated zeros in $B_\delta$. Let $m_i$ be
the number of them counted
with multiplicities. Then ${\ind}^{(i)}\Omega=m_i$, ${\ind}_P\Omega=
m_0-m_1+\ldots+(-1)^\ell m_\ell$.

{\bf 2.} Let $K_i$ ($i=0, 1, \ldots, \ell$) be the link
of the {\icis} $V\cap D_1\cap\ldots\cap D_i$, i.e., the intersection
$V\cap D_1\cap\ldots\cap D_i\cap S_\delta$. Let $d_i$ be the degree
of the map
$$
(\omega, df_1,\ldots, df_k,\ldots, df_{k+i}):K_i\to W_{k+i+1}(\CC^{n+k})
$$
($W_{k+i+1}(\CC^{n+k})$ is the Stiefel manifold of $(k+i+1)$--frames
in the dual $\CC^{n+k}$). Then ${\ind}^{(i)}\Omega=d_i$, ${\ind}_P\Omega=
d_0-d_1+\ldots+(-1)^\ell d_\ell$.

{\bf 3.} Let $I_i$ ($i=0, 1, \ldots, \ell$) be the ideal
of the ring ${\mathcal O}_{\CC^{n+k},0}$ of germs of holomorphic functions
of $n+k$ variables at the origin generated by $f_1$, \dots, $f_k$, \dots,
$f_{k+i}$ and the $(k+i+1)\times(k+i+1)$--minors of the matrix
$$\left( \begin{array}{ccc} \frac{\partial f_1}{\partial x_1} &
\cdots &
\frac{\partial f_1}{\partial x_{n+k}} \\
\vdots & \ddots & \vdots \\
\frac{\partial f_{k+i}}{\partial x_1} & \cdots & \frac{\partial
f_{k+i}}{\partial
x_{n+k}}\\
A_1 & \cdots & A_{n+k}
\end{array} \right).$$
Let $\nu_i =\dim_\CC{\mathcal O}_{\CC^{n+k},0}/I_i$. Then ${\ind}^{(i)}\Omega =
\nu_i$,
${\ind}_P\omega = \nu_0 -\nu_1 + \ldots + (-1)^\ell \nu_\ell$.

One has the following Poincar\'e--Hopf type formula.

\begin{theorem}\label{sec7theo2}
$$
\sum_P {\ind}_P\Omega=(-1)^n\chi(\widetilde V)
$$
where $\widetilde V$ is a smoothing of the variety $V$.
\end{theorem}

In order to get the Euler characteristic of the variety $V$ itself one
has to correct this formula by taking the Milnor numbers of the singular
points of $V$ into account. For $P \in V$, let $\mu_P$ denote the Milnor
number of the {\icis} $V$ at the point $P$ (see, e.g., \cite{Looijenga84}).
Note that, if $V$ is non-singular at $P$, then $\mu_P=0$.

\begin{theorem} \label{sec7theo2'}
$$
(-1)^n \chi(V) = \sum_P ({\ind}_P \Omega - \mu_P).
$$
\end{theorem}

\begin{remark} The Milnor number $\mu_P$ can also be written as an
alternating sum of indices of holomorphic 1-forms on several {\icis}.
One can say that the term $\mu_P$ corresponds to the difference between
two possible definitions of the index as in \cite[Proposition~1.4]{SS98}.
\end{remark}


In some cases it is natural to consider the situation when the pole locus
of the meromorphic 1-form $\omega$ is a multiple of an irreducible divisor.
Examples are:
\begin{itemize}
\item[1)] A polynomial 1-form on $\CC^n$ is a meromorphic 1-form on
$\CC\PP^n$ the pole locus of which is a multiple of the infinite hyperplane
$\CC\PP_\infty^{n-1}$.
\item[2)] Let $f$ be a meromorphic function with the pole locus $kD$
where $D$ is an irreducible divisor (i.e., in a neighbourhood of each
point $f$ can be written as $\widetilde{f}/F^k$ where $\widetilde{f}$
is holomorphic and $F=0$ is a local equation of $D$). Then its
differential $df$ is a meromorphic 1-form with the pole locus $(k+1)D$
(i.e., in a neighbourhood of each point $df$ can be written as
$\widetilde{\omega}/F^{k+1}$ where $\widetilde{\omega}$ is a holomorphic
1-form). A polynomial function of degree $k$ on $\CC^n$ is a meromorphic
function on $\CC\PP^n$ with the pole locus $k\CC\PP^{n-1}_\infty$.
\end{itemize}

Let us give a version of Theorem~\ref{sec7theo1} for this case. The fact
that the pole divisor of a 1-form is multiple means that the corresponding
line bundle is a power of another one. Let $\omega$ be a meromorphic 1-form
on a compact complex manifold $M^n$, that is, a holomorphic section of
the bundle $T^\ast M \otimes L$, where $L=\lambda^k$, $k>1$. Suppose that
$D_1=D$, $D_2$, \dots , $D_n$ are zero divisors of holomorphic sections of
the line bundle $\lambda$ such that, for each $i=1, \ldots , n$,
$D_1 \cap \ldots \cap D_i$ is non-singular and the form $\omega$ itself
and its restrictions to $D_1 \cap \ldots \cap D_i$ have isolated zeros.
Let $m_0$ (respectively $m_i$, $i=1, \ldots, n$) be the number of zeros
of the 1-form $\omega$ (respectively of its restriction to
$D_1 \cap \ldots \cap D_i$) counted with multiplicities.

\begin{theorem} \label{sec7theo3}
$$
\langle c_n(T^\ast M),[M]\rangle=
(-1)^n\chi(M)= m_0 -km_1+ \ldots + (-1)^nkm_n\,.
$$
\end{theorem}

\begin{remark} In \cite[Theorem~3.1]{HS98} there is given a formula for
the sum of the residues corresponding to the singular points of the
foliation on the projective plane $\CC\PP^2$ given by $df=0$, where
$f$ is a polynomial function of degree $k$ on $\CC^2$ (which defines a
meromorphic function on $\CC\PP^2$). In our terms this is the formula for
the number $m_0$ of zeros of the meromorphic 1-form $df$, the pole locus
of which is $k+1$ times the infinite line. Thus Theorem~\ref{sec7theo3}
can be considered as a generalization of \cite[Theorem~3.1]{HS98} to higher
dimensions and to meromorphic 1-forms which are not, in general,
differentials of functions.
\end{remark}


\section{Indices of collections of 1-forms}\label{sec8} 

One can say that all the indices discussed above are connected with
the Euler characteristic (e.g. through the Poincar\'{e}--Hopf type
theorem). If $M$ is a complex analytic manifold of dimension $n$,
then its Euler characteristic $\chi(M)$ is the characteristic number
$\langle c_n(TM), [M]\rangle=(-1)^n\langle c_n(T^\ast M), [M]\rangle$,
where $TM$ is the tangent bundle of the manifold $M$, $T^\ast M$ is
the dual bundle, and $c_n$ is the corresponding Chern class. One can
try to find generalizations of some notions and statements about indices
of vector fields or of 1-forms for other characteristic numbers (different
from $\langle c_n(TM), [M]\rangle$ or $\langle c_n(T^\ast M), [M]\rangle$).
This was made in \cite{EGBLMS} and \cite{EG05b} for the GSV index and
for the Euler obstruction.

The top Chern class of a vector bundle is the (first) obstruction to
existence of a non-vanishing section. Other Chern classes are
obstructions to existence of linear independent collections of
sections. Therefore instead of 1-forms on a complex variety we
consider collections of 1-forms.

Let $\pi: E \to M$ be a complex
analytic vector bundle of rank $m$ over a complex analytic
manifold $M$ of dimension $n$. It is known that the
($2(n-k)$-dimensional) cycle Poincar\'e dual to the characteristic
class $c_k(E)$ ($k=1, \ldots , m$) is represented by the set of
points of the manifold $M$ where $m-k+1$ generic sections of the
vector bundle $E$ are linearly dependent (cf., e.g., \cite[p.~413]{GH78}).

Let ${\mathcal M}(p,q)$, $D_{p,q}$, and $W_{p,q}$ be the spaces defined in
Section~\ref{sec3}.
Let ${\bf k}=(k_1, \ldots , k_s)$ be a sequence of positive integers with
$\sum_{i=1}^s k_i = k$. Consider the space ${\mathcal M}_{m, {\bf k}}=
\prod_{i=1}^s {\mathcal M}(m,m-k_i+1)$ and the subvariety $D_{m, {\bf k}}=
\prod_{i=1}^s D_{m,m-k_i+1}$ in it. The variety $D_{m, {\bf k}}$ consists
of sets $\{A_i\}$ of $m \times (m-k_i+1)$ matrices such that
$\mbox{rk}\, A_i < m-k_i+1$ for each $i=1, \ldots , s$. Since
$D_{m, {\bf k}}$ is irreducible of codimension $k$, its complement
$W_{m, {\bf k}}= {\mathcal M}_{m, {\bf k}} \setminus D_{m, {\bf k}}$
is $(2k-2)$-connected, $H_{2k-1}(W_{m, {\bf k}}) \cong \ZZ$, and there
is a natural choice of a generator of the latter group. This choice
defines a degree (an integer) of a map from an oriented manifold of
dimension $2k-1$ to the manifold $W_{m, {\bf k}}$.

Let $(V, 0) \subset (\CC^N,0)$ be an $n$-dimensional isolated
complete intersection singularity ({\icis}) defined by equations
$f_1=\ldots=f_{N-n}=0$ ($f_i\in{\mathcal O}_{\CC^N,0}$). Let $f$ be the
analytic map
$(f_1, \ldots , f_{N-n}): (\CC^N,0) \to (\CC^{N-n},0)$ ($V=f^{-1}(0)$).
Let $\{ X_j^{(i)} \}$ be a collection of vector
fields on a neighbourhood of the origin in $(\CC^N,0)$ ($i=1, \ldots , s$;
$j=1, \ldots , n-k_i+1$; $\sum k_i = n$) which are tangent to the
{\icis} $(V,0)=\{f_1= \cdots = f_{N-n}=0\} \subset (\CC^N,0)$ at non-singular
points of $V$.
We say that a point
$p\in V\setminus \{0\}$ is non-singular for the collection $\{ X_j^{(i)} \}$
 on $V$ if at least for some $i$ the vectors $X_1^{(i)}(p), \ldots , X_{n-k_i+1}^{(i)}(p)$ are
linearly independent. Suppose that the collection $\{ X^{(i)}_j\}$ has
no singular points on $V$ outside of the origin in a neighbourhood of it.
Let $U$ be a neighbourhood of the origin in $\CC^N$ where all the functions
$f_r$ ($r=1, \ldots , N-n$) and the vector fields $X_j^{(i)}$ are defined
and such that the collection $\{ X_j^{(i)} \}$ has no singular
points on $(V \cap U) \setminus \{ 0\}$. Let
$S_\delta \subset U$ be a sufficiently small sphere around the origin
which intersects $V$ transversally and denote by $K=V \cap S_\delta$
the link of the {\icis} $(V,0)$. The manifold $K$ has a natural orientation
as the boundary of a complex analytic manifold. Let $\Psi_V$ be the mapping
from $V \cap U$ to ${\mathcal M}_{n,{\bf k}}$ which sends a point
$x \in V \cap U$ to the collection of $N \times (N-k_i+1)$-matrices
$$
\{ ( \mbox{grad}\, f_1(x), \ldots , \mbox{grad}\, f_{N-n}(x),
X_1^{(i)}(x), \ldots , X_{n-k_i+1}^{(i)}(x))\}, \quad i=1,\ldots, s.
$$
Here
$\mbox{grad}\, f_r$ is the gradient vector field of $f_r$ defined in
Section~\ref{sec3}. Its restriction $\psi_V$ to the link $K$ maps $K$ to the subset
$W_{N,{\bf k}}$.

\begin{definition}
The
index
${\ind}_{V,0}\{ X^{(i)}_j\}$ of the
collection of vector fields $\{ \omega^{(i)}_j \}$ on the {\icis} $V$
is the degree of the mapping $\psi_V : K \to W_{N, {\bf k}}$.
\end{definition}

For $s=1$, $k_1=n$, this index is the GSV index of a vector field on an {\icis} (see Section~\ref{sec3}).

Let $V \subset \CC\PP^N$ be an $n$-dimensional complete intersection with
isolated singular points, $V=\{f_1=\ldots=f_{N-n}=0\}$ where $f_i$ are
homogeneous functions in $(N+1)$ variables. Let $\{ X^{(i)}_j \}$ be
a collection of continuous vector fields on $\CC\PP^N$ which are tangent
to $V$. Let $\widetilde V$ be a smoothing of the complete intersection $V$,
i.e. $\widetilde V$ is defined by $N-n$ equations $\widetilde f_1=\ldots=
\widetilde f_{N-n}=0$ where the homogeneous functions $\widetilde f_i$
are small perturbations of the functions $f_i$ and $\widetilde V$
is smooth.

The usual description of Chern classes of a vector bundle as obstructions
to existence of several linear independent sections of the bundle implies
the following statement.

\begin{theorem} \label{sec8theo1X}
One has
$$
\sum_{p \in V} {\ind}_p \{X^{(i)}_j\} = \langle \prod_{i=1}^s
c_{k_i}(T \widetilde{V}), [ \widetilde{V} ] \rangle,
$$
where $\widetilde{V}$ is a smoothing of the complete intersection $V$.
\end{theorem}

Now let
$\{ \omega^{(i)}_j\}$ be a collection of (continuous) 1-forms on a
neighbourhood of the origin in $(\CC^N,0)$ with $i=1, \ldots, s$,
$j=1, \ldots , n-k_i+1$, $\sum k_i = n$. We say that a point
$p\in V\setminus \{0\}$ is non-singular for the collection
$\{\omega^{(i)}_j\}$ on $V$ if at least for some $i$ the restrictions
of the 1-forms $\omega^{(i)}_j(p)$ to the tangent space $T_p V$ are
linearly independent. Suppose that the collection $\{\omega^{(i)}_j\}$ has
no singular points on $V$ outside of the origin in a neighbourhood of it.
As above let $K=V\cap S_\delta$ be the link of the {\icis} $(V,0)$
(all the functions $f_r$ and the 1-forms $\omega^{(i)}_j$ are defined
in a neighbourhood of the ball $B_\delta$.
Let $\Psi_V$ be the mapping
from $V \cap U$ to ${\mathcal M}_{n,{\bf k}}$ which sends a point
$x \in V \cap U$ to the collection of $N \times (N-k_i+1)$-matrices
$$
\{ (df_1(x), \ldots , df_{N-n}(x), \omega_1^{(i)}(x), \ldots ,
\omega_{n -k_i+1}^{(i)}(x)) \}, \quad i=1, \ldots, s.
$$
Its restriction $\psi_V$ to the link $K$ maps $K$ to the subset
$W_{N,{\bf k}}$.

\begin{definition}
The {\em index} ${\ind}_{V,0}\{\omega^{(i)}_j\}$ of the
collection of 1-forms $\{ \omega^{(i)}_j \}$ on the {\icis} $V$
is the degree of the mapping $\psi_V : K \to W_{N, {\bf k}}$.
\end{definition}

One can easily see that the index ${\ind}_{V,0}\{\omega^{(i)}_j\}$
is equal to the intersection number of the germ of the image of the
mapping $\Psi_V$ with the variety $D_{N, {\bf k}}$. If all the 1-forms
$\omega^{(i)}_j$ are complex analytic, the mapping $\Psi_V$ is complex
analytic as well.

For $s=1$, $k_1=n$, this index is the GSV index of a 1-form
(Section~\ref{sec3}).

Let $V \subset \CC\PP^N$ be an $n$-dimensional complete intersection with
isolated singular points, $V=\{f_1=\ldots=f_{N-n}=0\}$ where $f_i$ are
homogeneous functions in $(N+1)$ variables. Let $L$ be a complex line
bundle on $V$ and let $\{ \omega^{(i)}_j \}$ be a collection of continuous
1-forms on $V$ with values in $L$. Here this means that the forms
$\omega_j^{(i)}$ are continuous sections of the vector bundle
$T^\ast V \otimes L$ outside of the singular points of $V$. Since, in
a neighbourhood of each point $p$, the vector bundle $L$ is trivial,
one can define the index ${\ind}_p \{\omega^{(i)}_j\}$ of the
collection of 1-forms $\{ \omega_j^{(i)}\}$ at the point $p$ just in
the same way as in the local setting above.
Let $\widetilde V$ be a smoothing of the complete intersection $V$.
One can consider $L$ as a line bundle
on the smoothing $\widetilde{V}$ of the complete intersection $V$ as well
(e.g., using the pull back along a projection of $\widetilde{V}$ to $V$;
$L$ is not, in general, complex analytic).
The collection $\{\omega^{(i)}_j \}$ of 1-forms can also be extended
to a neighbourhood of $V$ in such a way that it will define a collection
of 1-forms on the smoothing $\widetilde{V}$ (also denoted by
$\{\omega^{(i)}_j \}$) with isolated singular points.
The sum of the indices of the collection $\{\omega^{(i)}_j \}$
on the smoothing $\widetilde{V}$ of $V$ in a neighbourhood of the point
$p$ is equal to the index ${\ind}_{V,p}\, \{\omega^{(i)}_j\}$.

One has the following analogue of Proposition \ref{prop_gsv3} for
1-forms.

\begin{theorem} \label{sec8theo1}
One has
$$
\sum_{p \in V} {\ind}_p \{\omega^{(i)}_j\} = \langle \prod_{i=1}^s
c_{k_i}(T^\ast\widetilde{V} \otimes L), [ \widetilde{V} ] \rangle,
$$
where $\widetilde{V}$ is a smoothing of the complete intersection $V$.
\end{theorem}

As above, let $(V,0) \subset (\CC^N,0)$ be the {\icis} defined by
the equations $f_1=\cdots = f_{N-n}=0$.
Let $\{ \omega_j^{(i)}\}$
($i=1, \ldots, s$; $j=1, \ldots , n-k_i+1$) be a collection of
1-forms on a neighbourhood of the origin in $\CC^N$ without singular
points on $V\setminus\{0\}$ in a neighbourhood of the origin. If all
the 1-forms $\omega_j^{(i)}$ are complex analytic, there exists an
algebraic formula for the index ${\ind}_{V,0}\{ \omega^{(i)}_j \}$
of the collection $\{ \omega^{(i)}_j \}$ similar to that from
Theorem~\ref{sec6theoAlg}.

Let $I_{V,\{ \omega_j^{(i)}\}}$ be
the ideal in the ring ${\mathcal O}_{\CC^N,0}$ generated by the functions
$f_1, \ldots , f_{N-n}$ and by the $(N-k_i+1) \times (N-k_i+1)$ minors
of all the matrices
$$
(df_1(x), \ldots , df_{N-n}(x), \omega_1^{(i)}(x), \ldots,
\omega_{n-k_i+1}^{(i)}(x))
$$
for all $i=1, \ldots, s$.

\begin{theorem} \label{sec8theo2} (see \cite{EGBLMS})
$$
{\ind}_{V,0}\{ \omega^{(i)}_j \} =
\dim_\CC {\mathcal O}_{\CC^N,0}/I_{V,\{\omega_j^{(i)}\}}.
$$
\end{theorem}

\begin{remark}
A formula similar to that of
Theorem~\ref{sec8theo2} does not exist for collections of vector fields.  A reason is that in this case
the index is the intersection number with $D_{N, {\bf k}}$ of the image
of the {\icis} $(V,0)$ under a map which is not complex analytic. Moreover,
in some cases this index can be negative (see e.g.
\cite[Proposition~2.2]{GSV91}).
\end{remark}

\begin{example}
Let $V$ be the (singular) quadric $x_1^2+x_2^2+x_3^2=0$ in the
projective space $\CC\PP^3$ with the coordinates
$(x_0:x_1:x_2:x_3)$. Let the 1-forms from the collection
$\{\{\omega^{(1)}_1, \omega^{(1)}_2\},\{\omega^{(2)}_1,
\omega^{(2)}_2\}\}$ be defined in the affine chart
$\CC^3=\{x_0=1\}$ of the projective space $\CC\PP^3$ as
$\omega^{(1)}_1=dx_1$, $\omega^{(1)}_2=x_2dx_3-x_3dx_2$,
$\omega^{(2)}_1=dx_2$, $\omega^{(2)}_2=x_1dx_3-x_3dx_1$. They are
holomorphic 1-forms on the projective space $\CC\PP^3$ with values
in the line bundle ${\mathcal O}(2)$. It is easy to see that being
restricted to the quadric $V$ the indicated collection of 1-forms
has no singular points outside of the singular point $(1:0:0:0)$
of $V$ itself. Theorem~\ref{sec8theo1} says that the index of this
point is equal to $\langle c_1^2\,(T^*\widetilde V\otimes
i^\ast{\mathcal O}(2)),[\widetilde{V}]\rangle$, where $\widetilde{V}$ is a
non-singular quadric, $i^\ast{\mathcal O}(2)$ is the restriction of the
line bundle ${\mathcal O}(2)$ to $\widetilde{V}$. Taking into account that
$\widetilde{V}\cong(\CC\PP^1)^2$, one gets that this
characteristic number is equal to $8$. Now by Theorem~\ref{sec8theo2}
the index of the indicated collection of 1-forms on the quadric
$V$ at the origin in $\CC^3\subset\CC\PP^3$ is equal to $\dim_\CC
{\mathcal O}_{\CC^3,0}/\langle x_1^2+x_2^2+x_3^2, x_2^2+x_3^2, x_1^2+x_3^2 \rangle$
which is also equal to $8$.
\end{example}


There exists a generalization of the notion of the Euler obstruction to
collections of 1-forms corresponding to different Chern numbers.

Let $(V^n,0)\subset(\CC^N,0)$ be the germ of a purely $n$-dimensional
reduced complex analytic variety at the origin (generally speaking with
a non-isolated singularity). Let ${\bf k}=\{k_i\}$, $i=1,\ldots, s$, be
a fixed partition of $n$ (i.e., $k_i$ are positive integers,
$\sum\limits_{i=1}^s k_i=n$). Let $\{\omega^{(i)}_j\}$ ($i=1,\ldots, s$,
$j=1,\ldots, n-k_i+1$) be a collection of germs of 1-forms on $(\CC^N, 0)$
(not necessarily complex analytic; it suffices that the forms
$\omega^{(i)}_j$ are complex linear functions on $\CC^N$ continuously
depending on
a point of $\CC^N$). Let $\eps>0$ be small enough so that there is
a representative $V$ of the germ $(V,0)$ and representatives
$\omega^{(i)}_j$ of the germs of 1-forms inside the ball
$B_\eps(0)\subset\CC^N$.

\begin{definition}
A point $P\in V$ is called a {\em special} point of the collection
$\{\omega^{(i)}_j\}$ of 1-forms on the variety $V$ if there exists
a sequence $\{P_m\}$ of points from the non-singular part $V_{\rm reg}$
of the variety $V$ such that the sequence $T_{P_m}V_{\rm reg}$ of the
tangent spaces at the points $P_m$ has a limit $L$ $($in $G(n,N)$$)$ and the
restrictions of the 1-forms $\omega^{(i)}_1$, \dots, $\omega^{(i)}_{n-k_i+1}$
to the subspace $L\subset T_P\CC^N$ are linearly dependent for each
$i=1, \ldots, s$. The collection $\{\omega^{(i)}_j\}$ of 1-forms has an
{\em isolated special point} on $(V,0)$ if it has no special points on
$V$ in a punctured neighbourhood of the origin.
\end{definition}

\begin{remark}
If the 1-forms $\omega^{(i)}_j$ are complex analytic, the property to have
an isolated special point is a condition on the classes of these 1-forms
in the module
$$
\Omega^1_{V,0}=
\Omega^1_{\CC^N,0}/\{f\cdot\Omega^1_{\CC^N,0} +
df\cdot{\mathcal O}_{\CC^N,0}\vert f\in{\mathcal J}_V\}
$$
of germs of 1-forms on the variety $V$ (${\mathcal J}_V$ is the ideal of
germs of functions vanishing on $V$).
\end{remark}

\begin{remark}
For the case $s=1$ (and therefore $k_1=n$), i.e. for one 1-form
$\omega$, we discussed the notion of a {\em singular} point of the 1-form
$\omega$ on $V$. One can easily see that a special point of the 1-form
$\omega$ on $V$ is singular, but not vice versa. (E.g.\ the origin is a
singular point of the 1-form $dx$ on the cone $\{x^2+y^2+z^2=0\}$, but
not a special one.) On a smooth variety these two notions coincide.
\end{remark}

\begin{definition}
A special (singular) point of a collection $\{\omega^{(i)}_j\}$ of
germs of 1-forms on a smooth $n$-dimensional variety $V$ is
{\em non-degenerate} if the map $\Psi_V: V \cap U \to {\mathcal M}_{n,{\bf k}}$
described above is transversal to
$D_{n, {\bf k}} \subset {\mathcal M}_{n, {\bf k}}$ at a non-singular point of
it.
\end{definition}

Let
$$
{\mathcal L}_{\bf k}= \prod\limits_{i=1}^s \prod\limits_{j=1}^{n-k_i+1}
\CC^{N\ast}_{ij}
$$
be the space of collections of linear functions on $\CC^N$ (i.e.\ of
1-forms with constant coefficients). The following statement holds.

\begin{proposition}\label{prop1a}
There exists an open and dense set $U \subset {\mathcal L}_{\bf k}$ such that
each collection $\{\ell^{(i)}_j\} \in U$ has only isolated special points
on $V$ and, moreover, all these points belong to the smooth part
$V_{\rm reg}$ of the variety $V$ and are non-degenerate.
\end{proposition}

\begin{corollary} \label{cor1}
Let $\{\omega^{(i)}_j\}$ be a collection of 1-forms on $V$ with an
isolated special point at the origin. Then there exists a deformation
$\{\widetilde \omega^{(i)}_j\}$ of the collection $\{\omega^{(i)}_j\}$
whose special points lie in $V_{\rm reg}$ and are non-degenerate. Moreover,
as such a deformation one can use $\{\omega^{(i)}_j+\lambda \ell^{(i)}_j\}$
with a generic collection $\{\ell^{(i)}_j\}\in {\mathcal L}_{\bf k}$.
\end{corollary}

\begin{corollary}
The set of collections of holomorphic 1-forms with a non-isolated
special point at the origin has infinite codimension in the space
of all holomorphic collections.
\end{corollary}

Let $\{\omega^{(i)}_j\}$ be a collection of germs of 1-forms on $(V, 0)$
with an isolated special point at the origin. Let $\nu : \widehat{V} \to V$
be the Nash transformation of the variety $V\subset B_\eps(0)$ (see
Section~\ref{sec5}). The collection of 1-forms $\{\omega^{(i)}_j\}$ gives
rise to a section $\widehat{\omega}$ of the bundle
$$
\widehat\TT=\bigoplus_{i=1}^s\bigoplus_{j=1}^{n-k_i+1}\widehat T^*_{i,j}
$$
where $\widehat{T}^\ast_{i,j}$ are copies of the dual Nash bundle
$\widehat{T}^\ast$ over the Nash transform $\widehat{V}$ numbered by
indices $i$ and $j$. Let $\widehat\DD\subset\widehat\TT$ be the set of
pairs $(x,\{\alpha^{(i)}_j\})$ where $x\in\widehat V$ and the collection
$\{\alpha^{(i)}_j\}$ of elements of $\widehat T_x^*$ (i.e.\ of linear
functions on $\widehat T_x$) is such that $\alpha^{(i)}_1$, \dots,
$\alpha^{(i)}_{n-k_i+1}$ are linearly dependent for each $i=1, \dots, s$.
The image of the section $\widehat\omega$ does not intersect $\widehat\DD$
outside of the preimage $\nu^{-1}(0)\subset\widehat V$ of the origin.
The map $\widehat\TT\setminus \widehat\DD\to \widehat V$ is a fibre
bundle. The fibre $W_x=\widehat\TT\setminus \widehat\DD$ of it is
$(2n-2)$-connected, its homology group $H_{2n-1}(W_x;\ZZ)$ is isomorphic
to $\ZZ$ and has a natural generator (see above). The latter fact implies
that the fibre bundle $\widehat\TT\setminus \widehat\DD\to \widehat V$
is homotopically simple in dimension $2n-1$, i.e.\ the fundamental group
$\pi_1(\widehat V)$ of the base acts trivially on the homotopy group
$\pi_{2n-1}(W_x)$ of the fibre, the last one being isomorphic to the
homology group $H_{2n-1}(W_x)$: see, e.g., \cite{Steenrod51}.

\begin{definition}
The {\em local Chern obstruction} ${\rm Ch}_{V,0}\,\{\omega^{(i)}_j\}$ of
the collections of germs of 1-forms $\{\omega^{(i)}_j\}$ on $(V,0)$ at the
origin is the (primary) obstruction to extend the section $\widehat{\omega}$
of the fibre bundle $\widehat\TT\setminus \widehat\DD\to \widehat V$ from
the preimage of a neighbourhood of the sphere $S_\eps= \partial B_\eps$ to
$\widehat V$, more precisely its value $($as an element of the homology group
$H^{2n}(\nu^{-1}(V\cap B_\eps), \nu^{-1}(V \cap S_\eps);\ZZ)$\,$)$ on the
fundamental class of the pair
$(\nu^{-1}(V\cap B_\eps), \nu^{-1}(V \cap S_\eps))$.
\end{definition}

The definition of the local Chern obstruction
${\rm Ch}_{V,0}\,\{\omega^{(i)}_j\}$ can be reformulated in the
following way. Let ${\mathcal D}^{\bf k}_V\subset \CC^N \times {\mathcal L}^{\bf k}$
be the closure of the set of pairs $(x, \{\ell^{(i)}_j\})$ such that
$x \in V_{\rm reg}$ and the restrictions of the linear functions
$\ell_1^{(i)}$, \dots , $\ell_{n-k_i+1}^{(i)}$ to
$T_x V_{\rm reg} \subset \CC^N$ are linearly dependent for each
$i=1, \ldots, s$. (For $s=1$, ${\bf k}=\{n\}$, ${\mathcal D}^{\bf k}_V$ is
the (non-projectivized) conormal space of $V$ \cite{Teissier82}.) The
collection $\{\omega^{(i)}_j\}$ of germs of 1-forms on $(\CC^N,0)$
defines a section $\check{\omega}$ of the {trivial} fibre bundle
$\CC^N \times {\mathcal L}^{\bf k} \to \CC^N$. Then
$$
{\rm Ch}_{V,0}\,\{\omega^{(i)}_j\} = ( \check{\omega}(\CC^N) \circ
{\mathcal D}^{\bf k}_V)_0
$$
where $(\cdot \circ \cdot )_0$ is the intersection number at the origin
in $\CC^N \times {\mathcal L}^{\bf k}$. This description can be considered as
a generalization of an expression of the local Euler obstruction as a
micro-local intersection number defined in \cite{KS90}, see also
\cite[Sections 5.0.3 and 5.2.1]{Schurmann03} and \cite{Schurmann04}.

\begin{remark}
On a smooth manifold $V$ the local Chern obstruction
${\rm Ch}_{V,0}\,\{\omega^{(i)}_j\}$ coincides with the index
${\ind}_{V,0} \, \{\omega^{(i)}_j\}$ of the collection
$\{\omega^{(i)}_j\}$ defined above.
\end{remark}

\begin{remark}
The local Euler obstruction is defined for vector fields as well
as for 1-forms. One can see that vector fields are not well adapted to
a definition of the local Chern obstruction. A more or less direct version
of the definition above for vector fields demands to consider vector
fields on a singular variety $V \subset \CC^N$ to be sections $X=X(x)$
of $T\CC^N_{|V}$ such that $X(x) \in T_x V \subset T_x\CC^N$ ($\dim T_x V$
is not constant). (Traditionally vector fields tangent to smooth strata
of the variety $V$ are considered.) There exist only continuous
(non-trivial, i.e.\ with $s>1$) collections of such vector fields ''on
$V$'' with isolated special points, but not holomorphic ones.
\end{remark}

\begin{remark}
The definition of the local Chern obstruction
${\rm Ch}_{V,0}\,\{\omega^{(i)}_j\}$ may also be formulated in terms
of a collection $\{\omega^{(i)}\}$ of germs of 1-forms with values in
vector spaces $L_i$ of dimensions $n-k_i+1$. Therefore (via differentials)
it is also defined for a collection $\{f^{(i)}\}$ of germs of maps
$f^{(i)}: (\CC^N,0) \to (\CC^{n-k_i+1},0)$ (just as the Euler
obstruction is defined for a germ of a function).
\end{remark}

Being a (primary) obstruction, the local Chern obstruction satisfies the
law of conservation of number, i.e.\ if a collection of 1-forms
$\{\widetilde\omega^{(i)}_j\}$ is a deformation of the collection
$\{\omega^{(i)}_j\}$ and has isolated special points on $V$, then
$$
{\rm Ch}_{V,0} \, \{\omega^{(i)}_j\} = \sum{\rm Ch}_{V,Q}\,
\{\widetilde\omega^{(i)}_j\}
$$
where the sum on the right hand side is over all special points $Q$ of
the collection $\{\widetilde\omega^{(i)}_j\}$ on $V$ in a neighbourhood
of the origin. With Corollary~\ref{cor1} this implies the following
statements.

\begin{proposition}\label{prop2}
The local Chern obstruction ${\rm Ch}_{V,0} \, \{\omega^{(i)}_j\}$ of a
collection $\{\omega^{(i)}_j\}$ of germs of holomorphic 1-forms is equal
to the number of special points on $V$ of a generic (holomorphic)
deformation of the collection.
\end{proposition}

This statement is an analogue of Proposition~2.3 in \cite{STV05b}.

\begin{proposition}\label{prop3}
If a collection $\{\omega^{(i)}_j\}$ of 1-forms on a compact (say,
projective) variety $V$ has only isolated special points, then the sum
of local Chern obstructions of the collection $\{\omega^{(i)}_j\}$ at
these points does not depend on the collection and therefore is an
invariant of the variety.
\end{proposition}

It is reasonable to consider this sum as ($(-1)^n$ times) the corresponding
Chern number of the singular variety $V$.

Let $(V,0)$ be an isolated complete intersection singularity.
The fact that both the Chern obstruction and the index of a collection
$\{\omega^{(i)}_j\}$ of 1-forms satisfy the law of conservation of number
and they coincide on a smooth manifold yields the following statement.

\begin{proposition}\label{prop4a}
For a collection $\{\omega^{(i)}_j\}$ on an isolated complete intersection
singularity $(V,0)$ the difference
$$
{\ind}_{V,0}\, \{\omega^{(i)}_j\} - {\rm Ch}_{V,0} \, \{\omega^{(i)}_j\}
$$
does not depend on the collection and therefore is an invariant of the germ
of the variety.
\end{proposition}

Since by Proposition~\ref{prop3} ${\rm Ch}_{V,0} \, \{\ell^{(i)}_j\}=0$
for a generic collection $\{\ell^{(i)}_j\}$ of linear functions on
$\CC^N$, one has the following statement.

\begin{corollary}
One has
$$
{\rm Ch}_{V,0} \, \{\omega^{(i)}_j\} =
{\ind}_{V,0}\, \{\omega^{(i)}_j\} - {\ind}_{V,0}\,\{\ell^{(i)}_j\}
$$
for a generic collection $\{\ell^{(i)}_j\}$ of linear functions on $\CC^N$.
\end{corollary}




\end{document}